\documentclass{article}[12pt]
\usepackage{mathrsfs}
\usepackage{amsmath}
\usepackage{amssymb, amsmath}
\usepackage{graphicx}
\usepackage{indentfirst}
\usepackage{bm}
\hoffset - 10 mm
\catcode`\@=11 \@addtoreset{equation}{section}

\catcode`\@=12
\usepackage{color}

\newtheorem{Theorem}{Theorem}[section]
\newtheorem{Proposition}{Proposition}[section]
\newtheorem{Lemma}{Lemma}[section]
\newtheorem{Corollary}{Corollary}[section]

\newtheorem{Remark}{Remark}[section]

\newcommand{\newcom}{\newcommand}
\newcommand{\bTheorem}[1]{
\begin{Theorem} \label{T#1} }
\newcommand{\eT}{\end{Theorem}}

\newcommand{\bProposition}[1]{
\begin{Proposition} \label{P#1}}
\newcommand{\eP}{\end{Proposition}}

\newcommand{\bLemma}[1]{
\begin{Lemma} \label{L#1} }
\newcommand{\eL}{\end{Lemma}}

\newcommand{\bCorollary}[1]{
\begin{Corollary} \label{C#1} }
\newcommand{\eC}{\end{Corollary}}

\newcommand{\proof}{{\bf Proof.}}
\newcommand{\beq}{\begin{equation}}
\newcommand{\eeq}{\end{equation}}
\newcom{\ben}{\begin{eqnarray}}
\newcom{\een}{\end{eqnarray}}
\newcom{\beno}{\begin{eqnarray*}}
\newcom{\eeno}{\end{eqnarray*}}
\newcom{\bali}{\begin{aligned}}
\newcom{\eali}{\end{aligned}}

\newcommand{\bFormula}[1]{
\begin{equation} \label{#1}}
\newcommand{\eF}{\end{equation}}

\newcommand{\p}{\partial}
\newcommand{\ls}{\lesssim}

\newcommand{\dx}{{\rm d} x}
\newcommand{\dy}{{\rm d} y}
\newcommand{\dt}{{\rm d} t }
\newcommand{\dtau}{{\rm d} \tau }
\newcommand{\dz}{{\rm d} z}
\newcommand{\dxi}{{\rm d}\xi}
\newcommand{\deta}{{\rm d}\eta}
\newcommand{\dzeta}{{\rm d}\zeta}
\newcommand{\dpsi}{{\rm d}\psi}
\newcommand{\dr}{{\rm d}r}
\newcommand{\dbeta}{{\rm d}\beta}
\newcommand{\dgamma}{{\rm d}\gamma}

\begin{document}

\title{\bf On asymptotic stability of the 3D Boussinesq equations without heat conduction}
\author{Lihua Dong \\ School of Physical and Mathematical Sciences,  Nanjing Tech University,\\ Nanjing 211816, China \\ Donglhmath@njtech.edu.cn \\ Yongzhong Sun \\ Department of Mathematics, Nanjing University, Nanjing 210093, China \\ sunyz@nju.edu.cn \\}
\maketitle

\abstract{We investigate the asymptotic stability of solution to Boussinesq equations without heat conduction with the initial data near a specific stationary solution in the three-dimensional domain $\Omega = \mathbb{R}^{2}\times (0,1)$. It is shown that the solution starting from a small perturbation to the stationary solution converges to it with explicit algebraic rates as time tends to infinity.}



\section{Introduction}

In the mathematical study of fluid dynamics, the Boussinesq approximation amounts to consider the motion of fluid  with density variations so small that it can be neglected in all terms except the one involving the gravity force. Note that the fluctuation of temperature is inversely proportional to that of the density. By using such an approximation, the Boussinesq system for incompressible and viscous fluid without heat conduction reads as follows.
\begin{equation}
     \label{Bou1}
     \left\{
     \begin{aligned}
       &\partial_{t}\mathbf{v} + \mathbf{v}\cdot \nabla \mathbf{v} - \Delta\mathbf{v} +\nabla p = \vartheta\mathbf{e}_{3},\\
      &\nabla\cdot\mathbf{v} = 0,\\
 &\partial_{t}\vartheta+\mathbf{v}\cdot \nabla \vartheta = 0.
     \end{aligned}
     \right.
\end{equation}
Here the unknowns $\mathbf{v}$, $p$ and $\vartheta$ are the velocity, pressure and temperature of the fluid, respectively. Note that both the viscosity coefficient and gravity constant are assumed to be $1$, which plays no special role in the following analysis. We refer to \cite{AJM1,AJM2,JP1} for detailed introduction on physical background of this and related models.

Due to rich physical background and mathematical features, in recent years there are many works on the Boussinesq equations without heat conduction from the point of view of rigorous mathematical analysis. In the two-dimensional case, the global existence, uniqueness and regularity of smooth solution have been investigated by many authors under different settings, see \cite{HA1,DC1,RD2,LH1,THSK1,TH1,WH1,NJ1,KW,ML1}. Moreover, there are some recent works on the asymptotic stability of certain stationary solution to (\ref{Bou1}), see \cite{CRD1,LD1,DS,LT1,RW2}. However, to the best of the authors' knowledge, there is no result on the asymptotic behavior in the three-dimensional case. In \cite{RD2}, R. Danchin and M. Paicu have obtained the global solution to (\ref{Bou1}) under the smallness assumption on the initial data. H. Abidi et al. \cite{HA2} and T. Hmidi et al. \cite{THR1} have shown global existence of large solution to (\ref{Bou1}) in the whole space $\mathbb{R}^3$ with axisymmetric initial data. Finally, we refer the reader to \cite{TE2,KW1} for works on the stability of certain stationary solution to the inviscid Boussinesq equations.

We assume that the fluid occupies the three-dimensional domain $\Omega=\mathbb{R}^{2}\times(0,1)\subset \mathbb{R}^{3}$ and
supplement system (\ref{Bou1}) with initial conditions
\begin{equation}\label{intiaL}
\mathbf{v}(0,\mathbf{x}) = \mathbf{v}_{0}(\mathbf{x}),\, \vartheta(0,\mathbf{x})=\vartheta_{0}(\mathbf{x})\text{ in }\Omega,
\end{equation}
together with the following slip-boundary conditions
\begin{equation}\label{boun1}
 (\mathbf{v}\cdot \mathbf{n})(t,\mathbf{x}) =0,\,\left(\mathbb{S}(\nabla\mathbf{v})\mathbf{n}\cdot\mathbf{\tau})(t,\mathbf{x}\right)=0 \text{ on }\partial\Omega,\,t>0.
\end{equation}
Here we use $\mathbf{x}=(x, y, z)$ to denote the spatial variable in $\mathbb{R}^{3}$, while $\mathbf{n}$  and $\tau$ are the unit outward normal and tangential direction to $\partial\Omega=\mathbb{R}^{2}\times\{z= 0,1\}$, respectively. The stress tensor $\mathbb{S}(\nabla\mathbf{v})=\frac{1}{2}\left(\nabla\mathbf{v}+\nabla^\top\mathbf{v}\right)$ is the symmetric part of $\nabla\mathbf{v}$. Since the boundary $\partial\Omega$ is flat, it follows that (\ref{boun1}) is equivalent to
\begin{equation}\label{boun21}
\partial_{3}v_{1}(t,\mathbf{x}) = 0,\, \partial_{3}v_{2}(t,\mathbf{x}) = 0,\,v_{3}(t,\mathbf{x})=0\text{ on }\partial\Omega,\,t>0.
\end{equation}

The purpose of this paper is to obtain the asymptotic stability of the specific stationary solution $(\mathbf{v}_s, p_{s}, \vartheta_{s})$ (see (\ref{ste1}) below) to (\ref{Bou1}) in the three-dimensional domain $\Omega=\mathbb{R}^{2}\times(0,1)$. We remark that if a stationary solution $\vartheta_{s}(z)$ satisfies $\vartheta'_{s}(z_0)<0$ for some $z_0\in [0,1]$, which implies that fluid with higher temperature lies below the one with lower temperature, then it is unstable--the Rayleigh-B\'enard instability happens, see \cite{Ch,PGD1}, among others. Hence, in the present work we focus on the opposite case $\vartheta'_s(z)>0$ for all $z\in [0,1]$, which implies that fluid of lower temperature lies below the fluid of higher temperature. Specifically, we choose
\begin{equation}\label{ste1}
\mathbf{v}_s = 0,\  p_{s} = \frac{1}{2}z^2,\, \vartheta_{s} = z,\  z \in [0,1].
\end{equation}
In the following Theorem \ref{Mrt1}, we show that as $t\to \infty$, the solution to (\ref{Bou1}) converges to the specific stationary solution (\ref{ste1}) with explicit rates in $t$ as long as the initial data is close enough to it in some suitable Sobolev spaces of high order.

By introducing the perturbation
\[
\mathbf{u} = \mathbf{v}-\mathbf{v}_s,\,q= p-p_s,\,\theta = \vartheta - \vartheta_{s},
\]
system (\ref{Bou1}) is transformed into
\begin{equation}\label{Bou2}
     \left\{
     \begin{aligned}
       &\partial_{t}\mathbf{u}+\mathbf{u}\cdot \nabla \mathbf{u} -\Delta\mathbf{u}  + \nabla q = \theta \mathbf{e}_3, \\
       &\nabla\cdot\mathbf{u} = 0,\\
       &\partial_{t}\theta+\mathbf{u}\cdot \nabla \theta = -u_{3}.
     \end{aligned}
     \right.
\end{equation}
The corresponding initial and boundary conditions for system (\ref{Bou2}) read as follows.
\begin{equation}\label{inibou1}
 \left\{
     \begin{aligned}
& \mathbf{u}(0,\mathbf{x}) = \mathbf{u}_{0}(\mathbf{x}),\, \theta(0,\mathbf{x})=\theta_{0}(\mathbf{x}) \text{ in }\Omega,\\
 &u_3(t,\mathbf{x}) =0,\  \p_3 u_1(t,\mathbf{x})=0,\,\p_3 u_2(t,\mathbf{x}) = 0 \text{ on }\partial\Omega,\,t>0.
 \end{aligned}
 \right.
\end{equation}
We eliminate the pressure $\nabla q$ by taking ${\rm curl} (=\nabla\times)$ on both side of (\ref{Bou2})$_1$ to reformulate it in terms of the vorticity $\boldsymbol{\omega} = \nabla \times \mathbf{u}$ as follows.
\begin{equation}\label{Bou3}
     \left\{
     \begin{aligned}
       &\partial_{t}\boldsymbol{\omega}-\Delta\boldsymbol{\omega} +\mathbf{u}\cdot \nabla \boldsymbol{\omega} - \boldsymbol{\omega}\cdot \nabla \mathbf{u} = \left(\partial_{2}\theta, -\partial_{1}\theta, 0\right),\\
       &\partial_{t}\theta+\mathbf{u}\cdot \nabla \theta = -u_{3},\\
       &\mathbf{u} = \nabla\times(-\Delta)^{-1}\boldsymbol{\omega},
     \end{aligned}
     \right.
\end{equation}
together with the initial and boundary conditions
\begin{equation}
     \label{inibou2}
     \left\{
     \begin{aligned}
       & (\boldsymbol{\omega}, \theta)(0,\mathbf{x}) = (\boldsymbol{\omega}_{0}, \theta_{0})(\mathbf{x}) \text{ in }\Omega,\\
       &  \omega_1(t,\mathbf{x}) = 0, \omega_2(t,\mathbf{x}) = 0, \p_3\omega_3(t,\mathbf{x}) = 0 \text{ on }\partial\Omega,\,t>0.
     \end{aligned}
     \right.
\end{equation}
Here we use $(-\Delta)^{-1}\bm\omega$ to denote the stream function $\bm\varphi=(\varphi_1,\varphi_2,\varphi_3)$ satisfying
\begin{equation}\label{vorphi1}
     \left\{
     \begin{aligned}
       & -\Delta\boldsymbol{\varphi} = \boldsymbol{\omega}\text{ in }\Omega,\\
       &  \varphi_{1}|_{\partial\Omega} = 0,\ \varphi_{2}|_{\partial\Omega} = 0,\ \partial_{3}\varphi_{3}|_{\partial\Omega} = 0.
     \end{aligned}
     \right.
\end{equation}
Also note that $\mathbf{u} =  \nabla \times \boldsymbol{\varphi}$ satisfies
\begin{equation}\label{vuw}
     \left\{
     \begin{aligned}
       &\nabla \times \mathbf{u} = \boldsymbol{\omega},\,\nabla\cdot\mathbf{u} = 0\text{ in }\Omega,\\
       & \partial_3{u}_1|_{\partial\Omega} = \partial_3{u}_2|_{\partial\Omega} = {u}_3|_{\partial\Omega} = 0.
     \end{aligned}
     \right.
\end{equation}

Before stating our main result we first give some notations used throughout this paper. For $m\in \mathbb{N}$ and $p\in [1,\infty]$, $W^{m,p}$ denotes the standard (inhomogeneous) Sobolev space on $\Omega$. Especially, we use $H^{m}$ to denote the $L^{2}$-based Sobolev spaces $W^{m,2}$ and  $\langle\cdot, \cdot \rangle$ to denote the inner product in $L^{2}$. Moreover, $\widehat{\Omega}$ is used to denote the frequency space $\mathbb{R}\times\mathbb{R}\times\mathbb{N}$ as well as $\widehat{L}^p=L^{p}(\widehat{\Omega})(1\leq p \leq \infty)$ the space of $p$-integrable functions $g(\xi,\eta,k)$, $(\xi,\eta,k)\in \widehat{\Omega}$, with the norm
\[
\|g\|_{\widehat{L}^p}:=\|g\|_{L^{p}(\widehat{\Omega})} =\left\{ \sum_{k\in \mathbb{N}}\int_{\mathbb{R}^{2}}|g(\xi,\eta,k)|^{p}\dxi\deta\right\}^{\frac{1}{p}}.
\]
In order to be compatible with the boundary conditions, it is necessary to introduce some additional spaces. To this end, we define for $m\in \mathbb{N}$,
\[
 \mathfrak{D}^{m} := \{f\in H^{m}: \partial_{3}^{n}f|_{\partial\Omega} = 0,\,n = 0, 2, \cdots, 2[(m-1)/2]\},
\]
\[
 \mathfrak{N}^{m} := \{f\in H^{m}: \partial_{3}^{n}f|_{\partial\Omega} = 0,\,n = 1, 3, \cdots, 2[m/2]-1\},
\]
\[
\mathcal{W}^{m} :=\{\mathbf{g}\in H^{m}\times H^{m}\times H^{m}: \mathbf{g}=(g_{1}, g_{2}, g_{3})\in \mathfrak{D}^{m}\times\mathfrak{D}^{m}\times\mathfrak{N}^{m}\},
\]
\[
\mathcal{V}^{m} :=\{\mathbf{g}\in H^{m}\times H^{m} \times H^{m}: \mathbf{g}=(g_{1}, g_{2}, g_{3})\in \mathfrak{N}^{m}\times\mathfrak{N}^{m}\times\mathfrak{D}^{m}\}.
\]
Finally, $A \lesssim B$ means that there exists a generic constant $C$ such that $A\leq CB$ while $A\sim B$ stands for $A\lesssim B$ and $B\lesssim A$.

Now we state the main result of this paper.
\begin{Theorem}\label{Mrt1}
Let $m\geq31$ be an integer. Assume that $\nabla\cdot \boldsymbol{\omega}_{0} = 0$ and
\[
\boldsymbol{\omega}_{0} =(\omega_{10}, \omega_{20}, \omega_{30}) \in \mathcal{W}^{m}\cap W^{7,1},
\]
\[
\Lambda^{-1}\omega_{30}\in L^2\cap L^1,\,\theta_{0}\in \mathfrak{D}^{m+1}\cap W^{10,1}.
\]
There exists $\epsilon_{0}>0$ depending only on $m$ such that if
\[
\|\theta_{0}\|_{W^{10,1}} + \|\theta_{0}\|_{H^{m+1}} +\|\boldsymbol{\omega}_{0}\|_{W^{7,1}} + \|\boldsymbol{\omega}_{0}\|_{H^{m}}+\|\Lambda^{-1}\omega_{30}\|_{L^{2}} + \|\Lambda^{-1}\omega_{30}\|_{L^{1}}\leq \epsilon_{0},
\]
then (\ref{Bou3})--(\ref{inibou2}) admits a unique global smooth solution
\[
\boldsymbol{\omega} \in C([0,\infty);\mathcal{W}^{m}),\,\theta \in C([0,\infty);\mathfrak{D}^{m+1})
\]
satisfying
\[
\|\boldsymbol{\omega}(t)\|_{H^m}^{2}+\|\theta(t)\|_{H^{m+1}}^{2} \ls \epsilon_{0}^{2},\text{ for all }t>0.
\]
Moreover,
\[
 \|\theta(t)\|_{H^{5}} \lesssim \langle t\rangle^{-\frac{1}{2}},\
\|\boldsymbol{\omega}(t)\|_{H^{3}}+\|\nabla_{h}\theta(t)\|_{H^{3}} \lesssim \langle t\rangle^{-1},
\]
\[
\|\nabla_{h}\boldsymbol{\omega}_h(t)\|_{H^{1}}+\|\nabla_{h}^{2}\theta(t)\|_{H^{1}} \lesssim \langle t\rangle^{-\frac{5}{4}},\|\omega_{3}(t)\|_{L^{\infty}} \lesssim \langle t\rangle^{-\frac{3}{2}},
\]
\[
\|\theta(t)\|_{L^{\infty}} + \|\p_{3}\theta(t)\|_{L^{\infty}}\lesssim \langle t\rangle^{-1},\,
\|\nabla_{h}\theta(t)\|_{L^{\infty}} + \|\boldsymbol{\omega}_h(t)\|_{L^{\infty}} \lesssim \langle t\rangle^{-\frac{5}{4}}.
\]
\end{Theorem}
\noindent Note that here and throughout this paper,  we denote $\langle t \rangle = \max\{1,t\}$ and the subscript $h$ of a vector field means its horizontal direction, namely, $\nabla_{h}=(\partial_{1}, \partial_{2})$, $\boldsymbol{\omega}_h=(\omega_{1}, \omega_{2})$ and $\mathbf{u}_{h}=(u_{1},u_{2})$, etc..
\begin{Remark}In fact, it holds
\[
\|\nabla_{h}\boldsymbol{\omega}_h(t)\|_{H^{1}}+\|\nabla_{h}^{2}\theta(t)\|_{H^{1}} \lesssim \langle t\rangle^{-\frac{3}{2}+},\,\|\nabla_{h}\theta(t)\|_{L^{\infty}} + \|\boldsymbol{\omega}_h(t)\|_{L^{\infty}} \lesssim \langle t\rangle^{-\frac{3}{2}+}.
\]
provided the initial data are smooth, that is, $m=\infty$. This will be obvious in the proof where interpolation inequalities are applied.
\end{Remark}
\begin{Remark}Going back to the original initial-boundary value problem (\ref{Bou1})--(\ref{boun1}), we have for all $t>0$,
\[
\|\vartheta(t)-z\|_{H^{5}} + \|\mathbf{v}_{h}(t)\|_{H^{4}} \lesssim \langle t\rangle^{-\frac{1}{2}},\,\|\nabla_{h}\left(\vartheta(t)-z\right)\|_{H^{3}}\ls \langle t\rangle^{-1},
\]
\[ \|\nabla_{h}^{2}\left(\vartheta(t)-z\right)\|_{H^{1}}+\|v_{3}(t)\|_{H^{3}}\lesssim \langle t\rangle^{-\frac{5}{4}},\
\|\nabla_{h}\left(\vartheta(t)-z\right)\|_{L^{\infty}} + \|\nabla\mathbf{v}(t)\|_{L^{\infty}}\lesssim \langle t\rangle^{-\frac{5}{4}},
\]
\[
\|\vartheta(t)-z\|_{L^{\infty}} + \|\p_{3}\left(\vartheta(t)-z\right)\|_{L^{\infty}} + \|\mathbf{v}(t)\|_{L^{\infty}}\lesssim \langle t\rangle^{-1}.
\]
Note that compared with the two-dimensional case in \cite{DS}, solutions to the 3D system (\ref{Bou1})--(\ref{boun1}) have better decay rates.
\end{Remark}

Since system (\ref{Bou1}) is absent of heat conduction, the action of buoyancy plays the key role in stabilization. As it is well-known, to obtain uniform estimates of solutions to (\ref{Bou3})--(\ref{inibou2}) in higher order Sobolev spaces, one has to show global-in-time $L^{1}$ estimate of $\|\nabla\mathbf{u}(t)\|_{L^{\infty}}$. It turns out that the crucial step is to give time decay estimates for $\theta$. By decoupling (\ref{Bou3}), we find that $\theta$ satisfies
\[
\partial_{tt}\theta - \Delta\partial_{t}\theta - \partial_{1}(-\Delta)^{-1}\partial_{1}\theta-\partial_{2}(-\Delta)^{-1}\partial_{2}\theta = f,
\]
where $f$ is a nonlinear term, see (\ref{if3}). Hence $\theta$ exhibits dissipation in two horizontal directions. The appearance of partial dissipation due to the action of buoyancy is a reminisce of early works (\cite{HAZ3,DengZ1}) on global well-posedness of MHD system in the absence of resistivity, where dissipation appears only in one direction-the direction of (strong) magnetic field. Although we have one additional direction of dissipation in the present case, the dissipation effect is weaker than that in MHD, which enforces us to introduce high order Sobolev spaces for the initial data. In fact, we obtain the time decay rate of $\|\partial_3\theta(t)\|_{L^{\infty}}$ is at most $\langle t\rangle^{-1}$ while $\|\nabla_h\theta(t)\|_{L^{\infty}}$ is at least $\langle t \rangle^{-\frac{5}{4}}$, which in turn implies that the decay rate of $\|\nabla\mathbf{u}(t)\|_{L^{\infty}}$ is $\langle t\rangle^{-\frac{5}{4}}$.  Similar to \cite{DS}, our method to prove Theorem \ref{Mrt1} is to show decay estimates for the linearized system by spectral analysis and then construct suitable energy functionals to the nonlinear system.

The rest of this paper is organized as follows. In the next section we give some preliminary results that will be used frequently. In Section \ref{ns} we prove the main result-Theorem \ref{Mrt1}-by first studying the linearized operator and then proving nonlinear stability with suitable choice of energy based on the linear analysis. Conclusions as well as discussion of possible future work are given in the final section.

\section{Preliminaries}\label{s2}

In this section, we give some necessary results that will be used later.

\subsection{Functional setting and Fourier expansion}

Motivated by \cite{AD1,AD2}, we assume the initial vorticity and temperature satisfy
\[
\omega_{j0}\in H^{m},\,\partial_{3}^{n}\omega_{j0}|_{\partial\Omega} = 0,\text{ for } j=1,2 \text{ and } n = 0, 2, \cdots, 2[(m-1)/2],
\]
\[
\omega_{30}\in H^{m},\ \partial_{3}^{n}\omega_{30}|_{\partial\Omega} = 0,
\text{ for } n = 1, 3, \cdots, 2[m/2]-1,
\]
\[
\theta_{0}\in H^{m+1},\ \partial_{3}^{n}\theta_{0}|_{\partial\Omega} = 0,
\text{ for } n = 0, 2, \cdots, 2[m/2].
\]
It turns out that these boundary conditions are propagated in time. In fact, from (\ref{Bou3})$_{2}$ and
$u_{3}(t)|_{\partial\Omega}=0$, we obtain the following transport equation on the boundary.
\[
\partial_{t}\theta(t)|_{\partial\Omega} + u_{1}(t)\partial_{1}\theta(t)|_{\partial\Omega} + u_{2}(t)\partial_{2}\theta(t)|_{\partial\Omega} = 0.
\]
Then one finds $\theta(t)|_{\partial\Omega} = 0$ by $\theta_{0}|_{\partial\Omega}=0$. Moreover, taking $\partial_{3}$ on $\nabla\cdot\mathbf{u} = 0$ gives rise to
\[
\partial_{3}^{2}u_{3}(t)|_{\partial\Omega} = -\partial_{3}\partial_{1}u_{1}(t)|_{\partial\Omega}-\partial_{3}\partial_{2}u_{2}(t)|_{\partial\Omega} = 0,
\]
where we use the fact that $\partial_{3}u_{1}(t)|_{\partial\Omega}=0$ and $\partial_{3}u_{2}(t)|_{\partial\Omega} = 0$. By using (\ref{Bou3})$_{1}$ and (\ref{inibou2})$_{2}$,
\[
\partial_{3}^{2}\omega_{1}(t)|_{\partial\Omega} = 0,\ \partial_{3}^{2}\omega_{2}(t)|_{\partial\Omega} = 0,\ \partial_{3}^{3}\omega_{3}(t)|_{\partial\Omega} = 0,
\]
which together with (\ref{Bou3})$_{3}$ implies that
 \[
\partial_{3}^{3}u_{1}(t)|_{\partial\Omega} = 0,~\partial_{3}^{3}u_{2}(t)|_{\partial\Omega} = 0,~\partial_{3}^{4}u_{3}(t)|_{\partial\Omega} = 0.
\]
Finally, taking $\partial_{3}^{2}$ on
(\ref{Bou3})$_2$ and restricting to the boundary yield
\begin{equation*}
  \partial_{t}\partial_{3}^{2}\theta(t)|_{\partial\Omega} + u_{1}(t)\partial_{1}\partial_{3}^{2}\theta(t)|_{\partial\Omega} +
  u_{2}(t)\partial_{2}\partial_{3}^{2}\theta(t)|_{\partial\Omega} + 2\partial_{3}u_{3}(t)\partial_{3}^{2}\theta(t)|_{\partial\Omega} = 0.
\end{equation*}
It follows that $\partial_{3}^{2}\theta(t)|_{\partial\Omega} = 0$ from $\partial_{3}^{2}\theta_{0}|_{\partial\Omega} = 0$. Similarly, from $\partial_{2}^{2n}\theta_{0}|_{\partial\Omega}=0$ for $n=2,4,\cdots$, we find that for sufficiently smooth solution $(\mathbf{u}, \theta)$,
\begin{equation}\label{P8}
  \begin{split}
& \partial_{3}^{2n}\theta(t)|_{\partial\Omega} = 0,~\partial_{3}^{2n}\omega_{1}(t)|_{\partial\Omega} = 0,\ \partial_{3}^{2n}\omega_{2}(t)|_{\partial\Omega} = 0,\ \partial_{3}^{2n+1}\omega_{3}(t)|_{\partial\Omega} = 0,  \\
&\partial_{3}^{2n+3}u_{1}(t)|_{\partial\Omega}= 0,\ \partial_{3}^{2n+3}u_{2}(t)|_{\partial\Omega}= 0,\ \partial_{3}^{2n +4}u_{3}(t)|_{\partial\Omega} = 0.
   \end{split}
\end{equation}
The argument above explains the choice of functional spaces for solutions $\mathbf{u}$ and $\theta$ in Theorem \ref{Mrt1}.

The Fourier expansion for a function $f\in \mathfrak{D}^m$ reads as
\begin{equation}\label{pr1}
f(x,y,z) = \frac{1}{\sqrt{2}\pi}\sum_{k=1}^{+\infty}\int_{\mathbb{R}^{2}}\widehat{f}_{o}(\xi,\eta,k)e^{i x \xi}e^{i y\eta}\dxi\deta \sin k\pi z,
\end{equation}
where
\[
\widehat{f}_{o}(\xi,\eta,k) = \frac{1}{\sqrt{2}\pi}\int_{\mathbb{R}^{2}}\int_{0}^{1}f(x,y,z)e^{- i x \xi} e^{-i y\eta}\sin k\pi z \dz\dx \dy, \text{ for } (\xi,\eta, k)\in \widehat{\Omega}.
\]
As for $f\in \mathfrak{N}^m$,
\begin{equation}\label{pr2}
f(x,y,z) = \frac{1}{\sqrt{2}\pi}\sum_{k=0}^{+\infty}\int_{\mathbb{R}^{2}}\widehat{f}_{e}(\xi,\eta,k)e^{ i x\xi}e^{i y\eta}\dxi \deta \cos k\pi z
\end{equation}
with
\[
\widehat{f}_{e}(\xi,\eta,k)= \frac{1}{\sqrt{2}\pi}\int_{\mathbb{R}^{2}}\int_{0}^{1}f(x,y,z)e^{- i x\xi}e^{- i y\eta}\cos k\pi z  \dz\dx \dy ,\ \text{ for } (\xi, \eta, k)\in \widehat{\Omega}.
\]
More specifically, since $f\in \mathfrak{D}^m $ (or $\mathfrak{N}^{m}$) implies $\partial_3 f \in \mathfrak{N}^{m-1}$ (or $ \mathfrak{D}^{m-1}$), it follows that
\[
\widehat{(\partial_3f)_e}(\xi,\eta,k) = -k \pi \widehat{f_o}(\xi,\eta,k)\,\left(\text{ or }\widehat{(\partial_3f)_o}(\xi,\eta,k) = k\pi \widehat{f_e}(\xi,\eta,k)\right).
\]
For notation convenience, we use $\widehat{f}$ (or $ \mathcal{F}(f))$ to denote $\widehat{f}_{o}$ or $\widehat{f}_{e}$ provided that $f \in \mathfrak{D}^m$ or $\mathfrak{N}^m$. Accordingly,
\begin{equation}\label{eLL2}
\|f\|_{H^m}\sim \|(1+|\cdot|^2)^{\frac{m}{2}}\widehat{f}(\cdot)\|_{\widehat{L}^2},\,
\|f\|_{\dot{H}^m}\sim \||\cdot|^{m}\widehat{f}(\cdot)\|_{\widehat{L}^2},\,\,f \in \mathfrak{D}^m \text{ or } \mathfrak{N}^m.
\end{equation}
In this paper, $\|(1+|\cdot|^2)^{\frac{m}{2}}\widehat{f}(\cdot)\|_{\widehat{L}^2}$ is denoted by $\|\widehat{f}(\cdot)\|_{\widehat{H}^m}$. Furthermore,
\begin{equation}\label{eLL3}
\|f\|_{L^\infty}\lesssim \|\widehat{f}\|_{\widehat{L}^1},\,\|\widehat{f}\|_{\widehat{L}^\infty}\lesssim \|{f}\|_{{L}^1},\,\|\widehat{f}\|_{\widehat{L}^1} \ls \|f\|_{H^2}\ls \|f\|_{W^{m,1}},\ m \ge 4,
\end{equation}
\begin{equation}\label{eLL4}
\|(1+|\cdot|^2)^{\frac{m}{2}}\widehat{f}\|_{\widehat{L}^\infty}\lesssim \|{f}\|_{{W}^{m,1}},\, f \in \mathfrak{D}^{m,1} \text{ or } \mathfrak{N}^{m,1}.
\end{equation}
Finally, $\Lambda^{\alpha}$ for $\alpha\in\mathbb{R}$ is defined as
\[
\widehat{\Lambda^{\alpha}f}(\xi,\eta,k) = \left(\xi^{2}+\eta^{2}+\pi^2k^{2}\right)^{\frac{\alpha}{2}}\widehat{f}(\xi,\eta,k),\,(\xi, \eta, k)\in \widehat{\Omega},\,f \in \mathfrak{D}^m \text{ or } \mathfrak{N}^m.
\]

\subsection{Basic lemmas}

We begin with some classical estimates on multiplication of two functions in Sobolev spaces, see (\cite{RA1,AJM2}), among others.
\begin{Lemma}\label{PL1}
 Let $m \in \mathbb{N}$.
\begin{itemize}
  \item If $f, g \in H^{m} \cap L^{\infty}$, then
\begin{equation}\label{BL11}
\|fg\|_{H^{m}}\lesssim \|f\|_{H^{m}}\|g\|_{L^{\infty}} + \|f\|_{L^{\infty}}\|g\|_{H^{m}}.
\end{equation}
\item
If $f \in H^{m} \cap W^{1,\infty}$, $g \in H^{m-1} \cap L^{\infty}$, then for $m\ge1$ and any $|\alpha|\leq m$,
\begin{equation}\label{BL14}
\|\partial^{\alpha}(fg) - f\partial^{\alpha}g\|_{L^{2}}\lesssim \|\nabla f\|_{L^{\infty}}\|g\|_{H^{m-1}} + \|f\|_{H^{m}}\|g\|_{L^{\infty}}.
\end{equation}
\end{itemize}
\end{Lemma}

The following elementary lemma is frequently used later, see \cite{RA1}, among others.
\begin{Lemma}\label{PL1m}
Let $m, m_1, m_2 \in \mathbb{N}$.

(i) Assume that $m_{1}\leq m \leq m_{2}$ and $f\in H^{m_{2}}$. Then
\begin{equation}\label{BLT1}
\|f\|_{H^{m}}\lesssim \|f\|_{H^{m_{1}}}^{s}\|f\|_{H^{m_{2}}}^{1-s},\text{ with } m= s m_{1} +(1-s)m_{2},\,0\leq s\leq1.
\end{equation}

(ii) If $f, g \in H^{m}$, then
\begin{equation}\label{BL13}
\|fg\|_{W^{m,1}}\lesssim \|f\|_{H^{m}}\|g\|_{L^{2}} + \|f\|_{L^{2}}\|g\|_{H^{m}}.
\end{equation}
\end{Lemma}

The next lemma follows from a direct calculation.
\begin{Lemma}\label{PL3}
Let $\mu, \nu>0$ be two constants such that $\mu \leq 1+\nu$. Then
\begin{equation}\label{BL22}
\int_{0}^{t}\frac{\dtau}{\langle t-\tau\rangle^{\mu}\langle \tau\rangle^{1+\nu}} \lesssim\langle t\rangle^{-\mu},\,
\int_{0}^{t}e^{-( t-\tau)}\langle \tau\rangle^{-\nu}\dtau \lesssim\langle t\rangle^{-\nu}.
\end{equation}
\end{Lemma}

Next we consider the following Laplace equation with Neumann boundary condition.
\begin{equation}
     \label{Po1}
     \left\{
     \begin{aligned}
       & -\Delta \psi = f \text{ in }\Omega,  \\
       &  \p_3\psi|_{\partial\Omega} = 0.
     \end{aligned}
     \right.
\end{equation}

\begin{Lemma}\label{PoL2}
Let $m \ge 0$ be an integer. Assume that
$f \in \mathfrak{N}^{m}$ and $\Lambda^{-1}f \in L^2$. Then there exists a unique solution (up to a constant) $\psi$ to (\ref{Po1}) such that
$\nabla\psi \in \mathcal{V}^{m+1}$.
Moreover,
\begin{equation}\label{Po3}
\|\nabla\psi\|_{H^{m+1}}\lesssim \|f\|_{H^{m}}+\|\Lambda^{-1}f \|_{L^{2}}.
\end{equation}
\end{Lemma}
{\bf Proof.} We search the solution in the form of
\[
\psi(x,y,z) = \frac{1}{\sqrt{2}\pi}\sum_{k=0}^{+\infty}\int_{\mathbb{R}^{2}}  \widehat{\psi}(\xi,\eta,k) e^{ i x\xi}e^{i y\eta}\dxi \deta \cos k\pi z.
\]
It is enough to take
\[  \widehat{\psi}(\xi,\eta,k) = \left({\xi^{2}+\eta^2+\pi^{2}k^{2}}\right)^{-1}\widehat{f}(\xi,\eta,k), \text{ for }(\xi,\eta, k)\in\widehat{\Omega}.
\]
Estimate (\ref{Po3}) follows from (\ref{eLL2}) directly.
\hfill$\square$

Similar result holds for the Dirichlet problem. More precisely, we have
\begin{Lemma}\label{PoL2D}
Let $m \ge 0$ be an integer and $f \in \mathfrak{D}^{m}$. Then there exists a unique solution $\psi\in\mathfrak{D}^{m+2}$ to
\begin{equation}\label{dd2}
     \left\{
     \begin{aligned}
       & -\Delta \psi = f \text{ in }\Omega,  \\
       &  \psi|_{\partial\Omega} = 0,
     \end{aligned}
     \right.
\end{equation}
such that
\begin{equation}\label{Po3D}
\|\psi\|_{H^{m+2}}\lesssim \|f\|_{H^{m}}.
\end{equation}
\end{Lemma}
For simplicity, we use $(-\Delta)^{-1}f$ to denote $(-\Delta_{N})^{-1}f$ or $(-\Delta_{D})^{-1}f$, which is the solution to the Neumann ((\ref{Po1})) or Dirichlet ((\ref{dd2})) problem, respectively.

Based on Lemma \ref{PoL2} and \ref{PoL2D}, we obtain from  (\ref{vorphi1}) the following
\begin{Lemma}\label{corollary1}
Let $m \in \mathbb{N}$. Given $\boldsymbol{\omega} \in \mathcal{W}^m$ such that $\Lambda^{-1}\omega_3\in L^2$, $\nabla\cdot\boldsymbol{\omega} = 0\text{ in }\Omega$, there exists a unique $\mathbf{u}\in \mathcal{V}^{m+1}$ satisfying (\ref{vuw}) such that
\begin{equation}\label{CE1s}
\|u_{3}\|_{H^{m+1}} + \|\nabla\mathbf{u}\|_{H^{m}}\lesssim \|\boldsymbol{\omega}\|_{H^{m}},\,
\end{equation}
\begin{equation}\label{CE1}
\|u_{1}\|_{H^{m+1}} + \|u_{2}\|_{H^{m+1}}\lesssim \|\boldsymbol{\omega}\|_{H^{m}} + \|\Lambda^{-1}\omega_3\|_{L^2}.
\end{equation}
\end{Lemma}

\section{Nonlinear stability}\label{ns}

The local existence result of (\ref{Bou3})--(\ref{inibou2}) can be proved by using the classical methods, see \cite{AD1,AD2,AJM2}. Here we omit the details for reasons of brevity.
\begin{Proposition}\label{Lth1}
Let $m \geq 3$ be an integer. Assume that $\boldsymbol{\omega}_{0}\in \mathcal{W}^{m}, \theta_{0}\in \mathfrak{D}^{m+1}$ and $\Lambda^{-1}\omega_{30}\in L^{2}$. There exists $T^{*}\in(0,+\infty]$ such that (\ref{Bou3})--(\ref{inibou2}) admits a unique solution
\[
(\boldsymbol{\omega}, \Lambda^{-1}\omega_{3}, \theta) \in C([0,T^{*}); \mathcal{W}^{m})\times C([0,T^{*}); L^{2})\times C([0,T^{*});\mathfrak{D}^{m+1}).
\]
\end{Proposition}
\begin{Remark}
From (\ref{P8}) and (\ref{CE1s})--(\ref{CE1}) in Lemma \ref{corollary1}, it follows that
\[
u_{1}\in C([0,T^{*});\mathfrak{N}^{m+1}),\,u_{2}\in C([0,T^{*});\mathfrak{N}^{m+1}),\,u_{3}\in C([0,T^{*});\mathfrak{D}^{m+1}).
\]
\end{Remark}

Based on Proposition \ref{Lth1}, we shall prove uniform-in-time estimates of the solution $(\boldsymbol{\omega}, \theta)$. First, from (\ref{Bou3})$_{1}$ we find that $\boldsymbol{\omega}$ satisfies the classical heat type equation with inhomogeneous terms. Now we decouple system (\ref{Bou3}) to get the equation for $\theta$. Taking the time derivative on (\ref{Bou3})$_{2}$ gives
\[
\partial_{tt} \theta + \p_tu_3 = - \p_t(\mathbf{u}\cdot\nabla\theta).
\]
Noting that
\[
\partial_t u_3 = \partial_{1}(-\Delta)^{-1}\partial_t\omega_{2} - \partial_{2}(-\Delta)^{-1}\partial_t\omega_{1},
\]
\[
-\Delta u_3 = \partial_t \Delta \theta + \Delta (\mathbf{u}\cdot\nabla\theta),
\]
we use (\ref{Bou3})$_{1}$ to find
\[
\partial_{tt}\theta-\Delta\partial_{t}\theta - \partial_{1}(-\Delta)^{-1}\partial_{1}\theta - \partial_{2}(-\Delta)^{-1}\partial_{2}\theta
\]
\[
=  -\partial_{t}\mathbf{u}\cdot\nabla\theta - \mathbf{u}\cdot\nabla\partial_{t}\theta + \partial_{1}(-\Delta)^{-1}(\mathbf{u}\cdot\nabla\omega_{2})
\]
\[
-\partial_{2}(-\Delta)^{-1}(\mathbf{u}\cdot\nabla\omega_{1})
+ \partial_{1}(-\Delta)^{-1}(\boldsymbol{\omega}\cdot\nabla u_{2})
\]
\beq\label{wave}
 - \partial_{2}(-\Delta)^{-1}(\boldsymbol{\omega}\cdot\nabla u_{1}) +
\Delta(\mathbf{u}\cdot\nabla\theta).
\eeq
To obtain decay estimates for $\theta$, we analyze the corresponding linearized equation in the following subsection.

\subsection{Analysis of the linearized equation for temperature}

To analyze (\ref{wave}), we consider the following inhomogeneous equation in $\Omega$,
\begin{equation}\label{Aa1}
   \partial_{tt}\phi-\Delta\partial_{t}\phi - \partial_{1}(-\Delta)^{-1}\partial_{1}\phi - \partial_{2}(-\Delta)^{-1}\partial_{2}\phi =F,
\end{equation}
together with the initial and boundary conditions
\begin{equation}\label{Aai1}
\left\{
     \begin{aligned}
       & \phi(0,\mathbf{x})=\phi_{0}(\mathbf{x}), \partial_{t}\phi(0,\mathbf{x})=\phi_{1}(\mathbf{x})\text{ in } \Omega ,\\
       &  \phi(t,\mathbf{x})=0 \text{ on } \partial\Omega,\,t>0.
     \end{aligned}
     \right.
\end{equation}
To solve it, we define operators $\mathcal{L}_{1}(t)$ and $\mathcal{L}_{2}(t)$ for $t>0$ as follows.
\begin{equation}\label{Aa2}
\mathcal{L}_{1}(t) = \frac{1}{2}\left(e^{-\frac{t}{2}\left(-\Delta-\chi\right)}+e^{-\frac{t}{2}\left(-\Delta +\chi\right)}\right),
\end{equation}
\begin{equation}\label{Aa3}
\mathcal{L}_{2}(t) = \frac{1}{\chi}\left(e^{-\frac{t}{2}\left(-\Delta -\chi\right)}-e^{-\frac{t}{2}\left(-\Delta+\chi\right)}\right).
\end{equation}
Here
\[
\chi=\sqrt{(-\Delta)^{2}+4(\partial_{1}(-\Delta)^{-1}\p_1 + \partial_{2}(-\Delta)^{-1}\p_2)}.
\]
\begin{Lemma}\label{AL1}
Assume that $\phi_0\in H^2\cap H_0^1$, $\phi_1\in L^2$ and $F\in L_{loc}^1(0, \infty; L^2)$. Then the solution to (\ref{Aa1})--(\ref{Aai1}) is given by
\begin{equation}\label{Aa4}
\phi(t,\mathbf{x})= \mathcal{L}_{1}(t)\phi_{0}(\mathbf{x}) + \mathcal{L}_{2}(t)\left(\frac{1}{2}(-\Delta)\phi_{0}(\mathbf{x}) +\phi_{1}(\mathbf{x})\right)
+ \int_{0}^{t}\mathcal{L}_{2}(t-\tau)F(\tau,\mathbf{x})\dtau.
\end{equation}
\end{Lemma}
\proof
Note that
\[
 \left[ \left(\partial_{t}- \frac{1}{2}\Delta\right)^{2}-\frac{1}{4}(-\Delta)^{2}-\partial_{1}(-\Delta)^{-1}\p_1 - \partial_{2}(-\Delta)^{-1}\p_2 \right]\phi(t,\mathbf{x})
\]
\[
=\left(\partial_{t}- \frac{1}{2}\Delta+\frac{1}{2}\sqrt{(-\Delta)^{2}+4(\partial_{1}(-\Delta)^{-1}\p_1 + \partial_{2}(-\Delta)^{-1}\p_2)}\right)
\]
\[
\cdot
\left(\partial_{t}- \frac{1}{2}\Delta-\frac{1}{2}\sqrt{(-\Delta)^{2}+4(\partial_{1}(-\Delta)^{-1}\p_1 + \partial_{2}(-\Delta)^{-1}\p_2)}\right)\phi(t,\mathbf{x}).
\]
Let
\begin{equation}\label{Aap1}
  \psi_{+}(t,\mathbf{x})=\left(\partial_{t}- \frac{1}{2}\Delta+\frac{1}{2}\chi\right)\phi(t,\mathbf{x})
\end{equation}
and
\begin{equation}\label{Aap2}
  \psi_{-}(t,\mathbf{x})=\left(\partial_{t}-\frac{1}{2}\Delta-\frac{1}{2}\chi\right)\phi(t,\mathbf{x}).
\end{equation}
Then
\[
\left(\partial_{t}- \frac{1}{2}\Delta-\frac{1}{2}\chi\right)\psi_{+}(t,\mathbf{x})=F(t,\mathbf{x}),\
\left(\partial_{t}- \frac{1}{2}\Delta+\frac{1}{2}\chi\right)\psi_{-}(t,\mathbf{x})=F(t,\mathbf{x}).
\]
From Duhamel's principle, it follows that
\beq\label{Aap3}
\psi_{+}(t,\mathbf{x}) = e^{t\left(\frac{1}{2}\Delta+\frac{1}{2}\chi\right)}\psi_{+}(0,\mathbf{x})
+\int_{0}^{t}e^{(t-\tau)\left(\frac{1}{2}\Delta+\frac{1}{2}\chi\right)}F(\tau,\mathbf{x})\dtau,
\eeq
\beq\label{Aap4}
\psi_{-}(t,\mathbf{x}) = e^{t\left(\frac{1}{2}\Delta-\frac{1}{2}\chi\right)}\psi_{-}(0,\mathbf{x})
+\int_{0}^{t}e^{(t-\tau)\left(\frac{1}{2}\Delta-\frac{1}{2}\chi\right)}F(\tau,\mathbf{x})\dtau.
\eeq
According to (\ref{Aap1}) and (\ref{Aap2}), we have
\beq\label{Aap5}
\psi_{+}(0,\mathbf{x})=\phi_{1}(\mathbf{x})-\frac{1}{2}\Delta\phi_{0}(\mathbf{x})+\frac{\chi}{2}\phi_{0}(\mathbf{x}),
\eeq
\beq\label{Aap6}
\psi_{-}(0,\mathbf{x})=\phi_{1}(\mathbf{x})-\frac{1}{2}\Delta\phi_{0}(\mathbf{x})-\frac{\chi}{2}\phi_{0}(\mathbf{x}).
\eeq
Substituting (\ref{Aap5}) and (\ref{Aap6}) into (\ref{Aap3}) and (\ref{Aap4}) respectively, then we use the fact that $\phi = \frac{\psi_{+}-\psi_{-}}{\chi}$ to obtain
\[
\phi(t,\mathbf{x}) = \frac{1}{2}\left(e^{-\frac{t}{2}\left(-\Delta+\chi\right)}+e^{-\frac{t}{2}\left(-\Delta-\chi\right)}\right)\phi_{0}(\mathbf{x})
\]
\[
+ \frac{1}{\chi}\left(e^{-\frac{t}{2}\left(-\Delta-\chi\right)}-e^{-\frac{t}{2}\left(-\Delta+\chi \right)}\right)\left(\phi_{1}(\mathbf{x})+\frac{1}{2}(-\Delta)\phi_{0}(\mathbf{x})\right)
\]
\[
 +\int_{0}^{t}\frac{1}{\chi}\left(e^{-\frac{t-\tau}{2}\left(-\Delta-\chi\right)}-e^{-\frac{t-\tau}{2}\left(-\Delta +\chi\right)}\right)F(\tau,\mathbf{x})\dtau.
\]
This completes the proof of Lemma \ref{AL1}. \hfill$\square$

The following lemma inspired by \cite{TE1} is used in the proof of decay estimates for $\mathcal{L}_{1}(t)$ and $\mathcal{L}_{2}(t)$.
\begin{Lemma}\label{spe2}
Let $g\in \mathfrak{D}^{8,1}$ and
\[
  \widehat{ G}(t,\xi,\eta,k) = e^{-\frac{\xi^{2}+\eta^2}{\left(\xi^{2}+\eta^2+\pi^{2}k^{2}\right)^{2}}t}\widehat{g}(\xi,\eta,k), \,(\xi,\eta,k) \in \widehat{\Omega}.
\]
Then
\begin{equation}\label{sph1}
\|\widehat{G}(t)\|_{\widehat{L}^{1}}\lesssim \langle t \rangle^{-1}\|g\|_{W^{6,1}},\  \|\widehat{\nabla_{h}G}(t)\|_{\widehat{L}^{1}}\lesssim \langle t \rangle^{-\frac{3}{2}}\|g\|_{W^{8,1}},
\end{equation}
\begin{equation}\label{sph2}
\|\widehat{G}(t)\|_{\widehat{L}^{2}}\lesssim \langle t \rangle^{-\frac{1}{2}}\|g\|_{W^{3,1}},\  \|\widehat{\nabla_{h}G}(t)\|_{\widehat{L}^{2}}\lesssim \langle t \rangle^{-1}\|g\|_{W^{5,1}}.
\end{equation}
\end{Lemma}
\proof
Note that
\[
\|\widehat{G}(t)\|_{\widehat{L}^{1}}\lesssim
\sum_{k=1}^{+\infty}\int_{\mathbb{R}^{2}}e^{-\frac{\xi^{2}+\eta^{2}}{\left(\xi^{2}+\eta^{2}+\pi^{2}k^{2}\right)^{2}}t}|\widehat{g}(\xi,\eta,k)|\dxi \deta
\]
\[
\ls
\|\left(\xi^{2}+\eta^{2}+\pi^{2}k^{2}\right)^{3}\widehat{g}(\xi,\eta,k)\|_{\widehat{L}^{\infty}}\int_{\pi}^{+\infty}\int_{\mathbb{R}^{2}}e^{-\frac{\xi^{2}+\eta^{2}}{\left(\xi^{2}+\eta^{2}+\zeta^{2}\right)^{2}}t}\left(\xi^{2}+\eta^{2}+\zeta^{2}\right)^{-3}\dxi \deta\dzeta
\]
\[
\lesssim \|g\|_{W^{6,1}}\int_{\pi}^{+\infty}\int_{\mathbb{R}^{2}}e^{-\frac{\xi^{2}+\eta^{2}}{\left(\xi^{2}+\eta^{2}+\zeta^{2}\right)^{2}}t}\left(\xi^{2}+\eta^{2}+\zeta^{2}\right)^{-3}\dxi\deta\dzeta.
\]
By applying spherical coordinates
\[
\xi = r\sin\beta\cos\psi,\ \pi\leq r<+\infty;\ \eta = r\sin\beta\sin\psi,\ 0\leq \beta\leq\pi;\ \zeta = r\cos\beta,\ 0\leq \psi\leq2\pi,
\]
then
\[
\int_{\pi}^{+\infty}\int_{\mathbb{R}^{2}}e^{-\frac{\xi^{2}+\eta^{2}}{\left(\xi^{2}+\eta^{2}+\zeta^{2}\right)^{2}}t}\left(\xi^{2}+\eta^{2}+\zeta^{2}\right)^{-3}\dxi\deta\dzeta
\lesssim \int_{\pi}^{+\infty}\int_{0}^{2\pi}\int_{0}^{\pi}e^{-\frac{\sin^{2}\beta}{ r^{2}}t}r^{-4}\sin\beta  \dbeta \dpsi \dr.
\]
We split the last integral as
\beq\label{deco}
\sum_{j=1}^{4}\int_{\pi}^{+\infty}\int_{0}^{2\pi}\int_{\frac{(j-1)\pi}{4}}^{\frac{j \pi}{4}}e^{-\frac{\sin^{2}\beta}{r^{2}}t}r^{-4}\sin\beta \dbeta \dpsi \dr = \sum_{j=1}^{4} Q_j.
\eeq
It is enough to consider the following two terms
\[
Q_{2}=\int_{\pi}^{+\infty}\int_{0}^{2\pi}\int_{\frac{\pi}{4}}^{\frac{\pi}{2}}e^{-\frac{\sin^{2}\beta}{ r^{2}}t}r^{-4}\sin\beta\dbeta \dpsi \dr,
\]
\[
Q_{4}=\int_{\pi}^{+\infty}\int_{0}^{2\pi}\int_{\frac{3\pi}{4}}^{\pi}e^{-\frac{\sin^{2}\beta}{ r^{2}}t}r^{-4}\sin\beta\dbeta \dpsi \dr.
\]
On one hand,
\[
Q_{2}\lesssim 1,\,Q_{4}\lesssim 1,\,
\frac{\sqrt{2}}{2}\leq |\sin \beta| \leq 1 \text{ for all }\beta \in \left[\frac{\pi}{4},\frac{\pi}{2}\right].
\]
On the other hand,
\[
Q_{2}
\lesssim \int_{\pi}^{+\infty}\int_{0}^{2\pi}\int_{\frac{\pi}{4}}^{\frac{\pi}{2}}e^{-\frac{t}{2r^{2}}}r^{-4}\dbeta \dpsi \dr
\]
\[
\lesssim \int_{\pi}^{+\infty}\int_{0}^{2\pi}\int_{\frac{\pi}{4}}^{\frac{\pi}{2}}e^{-\frac{t}{2r^{2}}}\frac{t}{2r^{2}}\frac{2r^{2}}{t}r^{-4} \dbeta \dpsi \dr
\lesssim t^{-1}.
\]
Hence $Q_2 \ls \langle t\rangle^{-1}$. Using change of variable $\gamma= \frac{\sqrt{t}}{\sqrt{2}r}\sin \beta$ gives
\[
Q_{4}
\lesssim t^{-\frac{1}{2}}\int_{\pi}^{+\infty}\int_{0}^{2\pi}\int_{\frac{3\pi}{4}}^{\pi}e^{-\frac{\sin^{2}\beta}{2 r^{2}}t}r^{-3}\dbeta \dpsi \dr
\]
\[\lesssim t^{-1}\int_{\pi}^{+\infty}\int_{0}^{\frac{\sqrt{t}}{2r}}e^{-\gamma^{2}}r^{-2}\frac{1}{\sqrt{1-\frac{2 r^{2}\gamma^{2}}{t}}}\dgamma \dr
\]
\[
\lesssim t^{-1}\int_{\pi}^{+\infty}\frac{1}{r^{2}}\dr\int_{-\infty}^{+\infty}e^{-\gamma^{2}}\dgamma
\lesssim t^{-1},
\]
where we use $e^{-\frac{\sin^{2}\beta}{2 r^{2}}t}\frac{\sin\beta}{r}t^{\frac{1}{2}}\ls 1$. Then $Q_4\ls \langle t \rangle^{-1}$. From (\ref{deco}), one infers that
\beq\label{r2d1}
\int_{\pi}^{+\infty}\int_{\mathbb{R}^{2}}e^{-\frac{\xi^{2}+\eta^{2}}{\left(\xi^{2}+\eta^{2}+\zeta^{2}\right)^{2}}t}\left(\xi^{2}+\eta^{2}+\zeta^{2}\right)^{-3}\dxi\deta\dzeta
\ls \langle t\rangle^{-1}.
\eeq
Thus
\[
\|\widehat{G}(t)\|_{\widehat{L}^{1}}\lesssim \langle t \rangle^{-1}\|g\|_{W^{6,1}}.
\]
Moreover,
\[
\|\widehat{\nabla_{h}G}(t)\|_{\widehat{L}^{1}}\lesssim
\sum_{k=1}^{+\infty}\int_{\mathbb{R}^{2}}\left(\xi^{2}+\eta^{2}\right)^{\frac{1}{2}}e^{-\frac{\xi^{2}+\eta^{2}}{\left(\xi^{2}+\eta^{2}+\pi^{2}k^{2}\right)^{2}}t}|\widehat{g}(\xi,\eta,k)|\dxi \deta
\]
\[
\ls
t^{-\frac{1}{2}}\sum_{k=1}^{+\infty}\int_{\mathbb{R}^{2}}t^{\frac{1}{2}}\frac{\left(\xi^{2}+\eta^{2}\right)^{\frac{1}{2}}}{\xi^{2}+\eta^{2}+\pi^{2}k^{2}}e^{-\frac{\xi^{2}+\eta^{2}}{\left(\xi^{2}+\eta^{2}+\pi^{2}k^{2}\right)^{2}}t}\left(\xi^{2}+\eta^{2}+\pi^{2}k^{2}\right)|\widehat{g}(\xi,\eta,k)|\dxi \deta
\]
\[
\ls
\|\left(\xi^{2}+\eta^{2}+\pi^{2}k^{2}\right)^{4}\widehat{g}(\xi,\eta,k)\|_{\widehat{L}^{\infty}}\int_{\pi}^{+\infty}\int_{\mathbb{R}^{2}}e^{-\frac{\xi^{2}+\eta^{2}}{\left(\xi^{2}+\eta^{2}+\zeta^{2}\right)^{2}}t}\left(\xi^{2}+\eta^{2}+\zeta^{2}\right)^{-3}\dxi \deta\dzeta
\]
\[
\lesssim t^{-\frac{1}{2}}\|g\|_{W^{8,1}}\int_{\pi}^{+\infty}\int_{\mathbb{R}^{2}}e^{-\frac{\xi^{2}+\eta^{2}}{\left(\xi^{2}+\eta^{2}+\zeta^{2}\right)^{2}}t}\left(\xi^{2}+\eta^{2}+\zeta^{2}\right)^{-3}\dxi\deta\dzeta.
\]
Using the fact that $\|\widehat{\nabla_{h}G}(t)\|_{\widehat{L}^{1}}\lesssim \|g\|_{W^{5,1}}$, we deduce from (\ref{r2d1}) that
\[
\|\widehat{\nabla_{h}G}(t)\|_{\widehat{L}^{1}}\lesssim \langle t \rangle^{-\frac{3}{2}}\|g\|_{W^{8,1}}.
\]
Note that
\[
\|\widehat{G}(t)\|_{\widehat{L}^{2}}^{2}\lesssim
\sum_{k=1}^{+\infty}\int_{\mathbb{R}^{2}}e^{-\frac{2(\xi^{2}+\eta^{2})}{\left(\xi^{2}+\eta^{2}+\pi^{2}k^{2}\right)^{2}}t}|\widehat{g}(\xi,\eta,k)|^{2}\dxi \deta
\]
\[
\ls
\|\left(\xi^{2}+\eta^{2}+\pi^{2}k^{2}\right)^{3}\widehat{g}^{2}(\xi,\eta,k)\|_{\widehat{L}^{\infty}}\int_{\pi}^{+\infty}\int_{\mathbb{R}^{2}}e^{-\frac{2(\xi^{2}+\eta^{2})}{\left(\xi^{2}+\eta^{2}+\zeta^{2}\right)^{2}}t}\left(\xi^{2}+\eta^{2}+\zeta^{2}\right)^{-3}\dxi \deta\dzeta
\]
\[
\lesssim \|g\|_{W^{3,1}}^{2}\int_{\pi}^{+\infty}\int_{\mathbb{R}^{2}}e^{-\frac{2(\xi^{2}+\eta^{2})}{\left(\xi^{2}+\eta^{2}+\zeta^{2}\right)^{2}}t}\left(\xi^{2}+\eta^{2}+\zeta^{2}\right)^{-3}\dxi\deta\dzeta.
\]
Hence
\[
\|\widehat{G}(t)\|_{\widehat{L}^{2}}\lesssim \langle t \rangle^{-\frac{1}{2}}\|g\|_{W^{3,1}}.
\]
Other estimate in (\ref{sph2}) can be obtained in the similar way.
\hfill$\square$

Based on Lemma \ref{spe2}, we now proceed to give decay estimates of $\mathcal{L}_{1}(t)$ and $\mathcal{L}_{2}(t)$.
\begin{Lemma}\label{AL2}
Let $f\in \mathfrak{D}^{8,1}$. Then
\begin{equation}\label{All1}
\left\|\mathcal{F}\left(\mathcal{L}_{1}(t)f\right)\right\|_{\widehat{L}^{1}} + \left\|\mathcal{F}\left(\mathcal{L}_{2}(t)f\right)\right\|_{\widehat{L}^{1}}\lesssim \langle t \rangle^{-1}\|f\|_{W^{6,1}},
\end{equation}
\begin{equation}\label{All1+}
\left\|\mathcal{F}\left(\partial_3\mathcal{L}_{1}(t)f\right)\right\|_{\widehat{L}^{1}} + \left\|\mathcal{F}\left(\partial_3\mathcal{L}_{2}(t)f\right)\right\|_{\widehat{L}^{1}}\lesssim \langle t \rangle^{-1}\|f\|_{W^{7,1}},
\end{equation}
\begin{equation}\label{All2}
\left\|\mathcal{F}\left(\partial_{t}\mathcal{L}_{1}(t)f\right)\right\|_{\widehat{L}^{1}} + \left\|\mathcal{F}\left(\partial_{t}\mathcal{L}_{2}(t)f\right)\right\|_{\widehat{L}^{1}}\lesssim \langle t \rangle^{-2}\|f\|_{W^{6,1}},
\end{equation}
\begin{equation}\label{All3}
\left\|\mathcal{F}\left(\nabla_{h}\mathcal{L}_{1}(t)f\right)\right\|_{\widehat{L}^{1}} + \left\|\mathcal{F}\left(\nabla_{h}\mathcal{L}_{2}(t)f\right)\right\|_{\widehat{L}^{1}}\lesssim \langle t \rangle^{-\frac{3}{2}}\|f\|_{W^{8,1}},
\end{equation}
\begin{equation}\label{All4}
\left\|\mathcal{F}\left(\mathcal{L}_{1}(t)f\right)\right\|_{\widehat{L}^{2}} + \left\|\mathcal{F}\left(\mathcal{L}_{2}(t)f\right)\right\|_{\widehat{L}^{2}}\lesssim \langle t \rangle^{-\frac{1}{2}}\|f\|_{W^{3,1}},
\end{equation}
\begin{equation}\label{All5}
\left\|\mathcal{F}\left(\partial_{t}\mathcal{L}_{1}(t)f\right)\right\|_{\widehat{L}^{2}} + \left\|\mathcal{F}\left(\partial_{t}\mathcal{L}_{2}(t)f\right)\right\|_{\widehat{L}^{2}}\lesssim \langle t \rangle^{-\frac{3}{2}}\|f\|_{W^{4,1}},
\end{equation}
\begin{equation}\label{All7}
\left\|\mathcal{F}\left(\nabla_{h}\mathcal{L}_{1}(t)f\right)\right\|_{\widehat{L}^{2}} + \left\|\mathcal{F}\left(\nabla_{h}\mathcal{L}_{2}(t)f\right)\right\|_{\widehat{L}^{2}}\lesssim \langle t \rangle^{-1}\|f\|_{W^{5,1}},
\end{equation}
\begin{equation}\label{All6}
\left\|\mathcal{F}\left(\nabla_{h}^{2}\mathcal{L}_{1}(t)f\right)\right\|_{\widehat{L}^{2}} + \left\|\mathcal{F}\left(\nabla_{h}^{2}\mathcal{L}_{2}(t)f\right)\right\|_{\widehat{L}^{2}}\lesssim \langle t \rangle^{-\frac{3}{2}}\|f\|_{W^{7,1}}.
\end{equation}
\end{Lemma}
\proof
Performing Fourier transform to (\ref{Aa2}) and (\ref{Aa3}),
\[
\widehat{\mathcal{L}_{1}(t)} = \frac{1}{2}\left(e^{-\frac{1}{2}\left(\xi^{2}+\eta^2+ \pi^{2}k^{2}+\sigma\right)t}+e^{-\frac{1}{2}\left(\xi^{2}+\eta^2+\pi^{2} k^{2}-\sigma\right)t}\right),
\]
\[
\widehat{\mathcal{L}_{2}(t)} = \frac{1}{2}\left(e^{-\frac{1}{2}\left(\xi^{2}+\eta^2+ \pi^{2}k^{2}-\sigma\right)t}-e^{-\frac{1}{2}\left(\xi^{2}+\eta^2+\pi^{2} k^{2}+\sigma\right)t}\right),
\]
where
\[
\sigma= \sqrt{(\xi^{2}+\eta^2+ \pi^{2}k^{2})^{2}-\frac{4(\xi^{2}+\eta^2)}{\xi^{2}+\eta^2+\pi^{2}k^{2}}}.
\]
Let
\[
\lambda_{+}= \frac{-(\xi^{2}+\eta^2+\pi^{2}k^{2})+\sigma}{2},\,\lambda_{-}= \frac{-(\xi^{2}+\eta^2+\pi^{2}k^{2})-\sigma}{2}.
\]
We calculate that
\[
\widehat{\mathcal{L}_{1}(t)} = \frac{1}{2}\left(e^{\lambda_{+}t}+e^{\lambda_{-}t}\right),\
\widehat{\mathcal{L}_{2}(t)}= \frac{1}{\lambda_{+}-\lambda_{-}}\left(e^{\lambda_{+}t}-e^{\lambda_{-}t}\right),
\]
\[
\widehat{\partial_{t}\mathcal{L}_{1}(t)} = \frac{1}{2}\left(\lambda_{+}e^{\lambda_{+}t}+\lambda_{-}e^{\lambda_{-}t}\right),\
\widehat{\partial_{t}\mathcal{L}_{2}(t)}= \frac{1}{\lambda_{+}-\lambda_{-}}\left(\lambda_{+}e^{\lambda_{+}t}-\lambda_{-}e^{\lambda_{-}t}\right),
\]
where $\lambda_{+}-\lambda_{-} = \sigma$.

Since $k\geq 1$, there exists $\lambda_{0}>0$ such that
$\lambda_{+} -\lambda_{-}\geq \lambda_{0}> 0$.
A straightforward calculation shows that
\[
\lambda_{+} = -\frac{1}{2}\left(\xi^{2}+\eta^2+\pi^{2}k^{2}-\sigma\right)
 \]
\[
= -\frac{1}{2}\left(\xi^{2}+\eta^2+\pi^{2}k^{2}-\sqrt{(\xi^{2}+\eta^2+\pi^{2}k^{2})^{2}-\frac{4(\xi^{2}+\eta^2)}{\xi^{2}+\eta^2+\pi^{2}k^{2}}}\right)
\]
\[
= -\frac{2(\xi^{2}+\eta^2)}{(\xi^{2}+\eta^2+\pi^{2}k^{2})^{2}}\frac{1}{1+\sqrt{1-\frac{4(\xi^{2}+\eta^2)}{(\xi^{2}+\eta^2+\pi^{2}k^{2})^{3}}}}.
\]
Moreover, we find
\[
-\frac{2(\xi^{2}+\eta^2)}{(\xi^{2}+\eta^2+\pi^{2}k^{2})^{2}}\leq \lambda_{+} \leq -\frac{\xi^{2}+\eta^2}{(\xi^{2}+\eta^2+\pi^{2}k^{2})^{2}},
\]
\[
-(\xi^{2}+\eta^2+\pi^{2}k^{2})\leq\lambda_{-} \leq -\frac{1}{2}(\xi^{2}+\eta^2+\pi^{2}k^{2}),
\]
\[
\lambda_{0}\leq\lambda_{+} -\lambda_{-} < \xi^{2}+\eta^2+\pi^{2}k^{2}.
\]
For $k\geq 1$,
\[
0<\widehat{\mathcal{L}_{1}(t)}\lesssim e^{-\frac{\xi^{2}+\eta^2}{(\xi^{2}+\eta^2+\pi^{2}k^{2})^{2}}t},
\]
\[
0\leq \widehat{\mathcal{L}_{2}(t)}\lesssim \frac{1}{\lambda_{0}}e^{-\frac{\xi^{2}+\eta^2}{(\xi^{2}+\eta^2+\pi^{2}k^{2})^{2}}t} + \frac{1}{\lambda_{0}}e^{-\frac{(\xi^{2}+\eta^2+\pi^{2}k^{2})}{2}t}
\lesssim e^{-\frac{\xi^{2}+\eta^2}{(\xi^{2}+\eta^2+\pi^{2}k^{2})^{2}}t},
\]
\[
\left|\widehat{\partial_{t}\mathcal{L}_{1}(t)}\right|\ls
\frac{\xi^{2}+\eta^2}{(\xi^{2}+\eta^2+\pi^{2}k^{2})^{2}}e^{-\frac{\xi^{2}+\eta^2}{(\xi^{2}+\eta^2+\pi^{2}k^{2})^{2}}t}
+ (\xi^{2}+\eta^2+\pi^{2}k^{2})e^{-\frac{\xi^{2}+\eta^2+\pi^{2}k^{2}}{2}t},
\]
\[
\left|\widehat{\partial_{t}\mathcal{L}_{2}(t)}\right|\ls
\frac{\xi^{2}+\eta^2}{(\xi^{2}+\eta^2+\pi^{2}k^{2})^{2}}e^{-\frac{\xi^{2} + \eta^2}{(\xi^{2}+\eta^2+\pi^{2}k^{2})^{2}}t}
+(\xi^{2}+\eta^2+\pi^{2}k^{2})e^{-\frac{\xi^{2}+\eta^2+\pi^{2}k^{2}}{2}t}.
\]

Note that
\[
\mathcal{F}\left(\mathcal{L}_{1}(t)f\right)(\xi,\eta,k) = \widehat{\mathcal{L}_{1}(t)}\widehat{f}(\xi,\eta,k)
\]
\[
 =
 \frac{1}{2}\left(e^{-\frac{1}{2}\left(\xi^{2}+\eta^2+ \pi^{2}k^{2}+\sigma\right)t}+e^{-\frac{1}{2}\left(\xi^{2}+\eta^2+\pi^{2} k^{2}-\sigma\right)t}\right)\widehat{f}(\xi,\eta,k).
\]
By using the fact that $0<\widehat{\mathcal{L}_{1}(t)}\lesssim e^{-\frac{\xi^{2}+\eta^2}{(\xi^{2}+\eta^2+\pi^{2}k^{2})^{2}}t}$ and applying Lemma \ref{spe2},
\[
\left\|\mathcal{F}\left(\mathcal{L}_{1}(t)f\right)\right\|_{\widehat{L}^{1}}
\lesssim
\sum_{k=1}^{\infty}\int_{\mathbb{R}^{2}}e^{-\frac{\xi^{2}+\eta^2}{(\xi^{2}+\eta^2+\pi^{2}k^{2})^{2}}t}|\widehat{f}(\xi,\eta,k)|\dxi\deta
\lesssim \langle t \rangle^{-1}\|f\|_{W^{6,1}}.
\]
Furthermore, we use Lemma \ref{spe2} to deduce that
\[
\left\|\mathcal{F}\left(\nabla_{h}\mathcal{L}_{1}(t)f\right)\right\|_{\widehat{L}^{1}}
\lesssim \sum_{k=1}^{+\infty}\int_{\mathbb{R}^{2}}e^{-\frac{\xi^{2} + \eta^2}{\left(\xi^{2}+\eta^2+ \pi^{2}k^{2}\right)^{2}}t}\left(\xi^{2}+\eta^2\right)^{\frac{1}{2}}|\widehat{f}(\xi,\eta,k)|\dxi\deta
\]
\[
\lesssim t^{-\frac{1}{2}}\sum_{k=1}^{+\infty}\int_{\mathbb{R}^{2}}e^{-\frac{\xi^{2} + \eta^2}{\left(\xi^{2}+\eta^2+ \pi^{2}k^{2}\right)^{2}}t}\frac{t^{\frac{1}{2}}\left(\xi^{2} + \eta^2\right)^{\frac{1}{2}}}{\xi^{2}+\eta^2+\pi^{2}k^{2}}\left(\xi^{2}+\eta^2+ \pi^{2}k^{2}\right)|\widehat{f}(\xi,\eta,k)|\dxi\deta
\]
\[
\lesssim t^{-\frac{1}{2}}\sum_{k=1}^{+\infty}\int_{\mathbb{R}^{2}}e^{-\frac{\xi^{2} + \eta^2}{2\left(\xi^{2}+\eta^2+ \pi^{2}k^{2}\right)^{2}}t}\left(\xi^{2}+\eta^2+ \pi^{2}k^{2}\right)|\widehat{f}(\xi,\eta,k)|\dxi\deta
 \lesssim \langle t \rangle^{-\frac{3}{2}}\|f\|_{W^{8,1}}.
\]
From (\ref{eLL3}) and $k\geq 1$, one infers that
\[
\left\|\mathcal{F}\left(\partial_{t}\mathcal{L}_{1}(t)f\right)\right\|_{\widehat{L}^{1}}
\ls \sum_{k=1}^{\infty}\int_{\mathbb{R}^{2}}\frac{\xi^{2} + \eta^2}{(\xi^{2}+\eta^2+\pi^{2}k^{2})^{2}}e^{-\frac{\xi^{2}+\eta^2}{\xi^{2}+\eta^2+\pi^2k^{2})^{2}}t}|\widehat{f}(\xi,\eta,k)|\dxi\deta
\]
\[
+ \sum_{k=1}^{\infty}\int_{\mathbb{R}^{2}}(\xi^{2}+\eta^2+\pi^{2}k^{2})e^{-\frac{\xi^{2}+\eta^2+\pi^{2}k^{2}}{2}t}|\widehat{f}(\xi,\eta,k)|\dxi\deta
\]
\[
\ls \langle t \rangle^{-2}\|f\|_{W^{6,1}} + e^{-t}\|(\xi^{2}+\eta^2+\pi^{2}k^{2})\widehat{f}\|_{\widehat{L}^{1}}
\ls \langle t \rangle^{-2}\|f\|_{W^{6,1}}.
\]
Other estimates in (\ref{All1})--(\ref{All6}) can be obtained in a similar way.
\hfill$\square$
\begin{Remark}
From Lemma \ref{AL2}, we observe that one derivative in the horizontal direction improves $\frac{1}{2}$-order of decay rate.
\end{Remark}
\subsection{Integral form of solutions}

Now we transform (\ref{Bou3}) into
\begin{equation}\label{if1}
     \left\{
     \begin{aligned}
       &\partial_{t}\boldsymbol{\omega}-\Delta\boldsymbol{\omega} = \mathbf{f}_{1},\\
        &\partial_{tt}\theta-\Delta\partial_{t}\theta - \partial_{1}(-\Delta)^{-1}\partial_{1}\theta - \partial_{2}(-\Delta)^{-1}\partial_{2}\theta = f_{2},\\
     \end{aligned}
     \right.
\end{equation}
with
\beq\label{if2}
\mathbf{f}_{1} = - \mathbf{u}\cdot \nabla \boldsymbol{\omega} + \boldsymbol{\omega}\cdot \nabla \mathbf{u} + \left(\partial_{2}\theta, -\partial_{1}\theta, 0\right),
\eeq
\[
  f_{2} =  -\partial_{t}\mathbf{u}\cdot\nabla\theta - \mathbf{u}\cdot\nabla\partial_{t}\theta +
\Delta(\mathbf{u}\cdot\nabla\theta)+\partial_{1}(-\Delta)^{-1}(\mathbf{u}\cdot\nabla\omega_{2})
\]
\beq\label{if3}
-\partial_{2}(-\Delta)^{-1}(\mathbf{u}\cdot\nabla\omega_{1}) + \partial_{1}(-\Delta)^{-1}(\boldsymbol{\omega}\cdot\nabla u_{2}) - \partial_{2}(-\Delta)^{-1}(\boldsymbol{\omega}\cdot\nabla u_{1}).
\eeq
From (\ref{P8}) and $\mathbf{f}_{1}= (f_{11},f_{12},f_{13})$, we check that
\[
f_{11} \in \mathfrak{D}^{m-1},\,f_{12} \in \mathfrak{D}^{m-1},\,f_{13} \in \mathfrak{N}^{m-1},\,f_2 \in \mathfrak{D}^{m-2}.
\]
Thus, by Duhamel's principle and Lemma \ref{AL1},
\begin{equation}\label{if4}
   \boldsymbol{\omega}(t,\mathbf{x}) = e^{t\Delta }\boldsymbol{\omega}_{0}(\mathbf{x}) + \int_{0}^{t}e^{(t-\tau)\Delta}\mathbf{f}_{1}(\tau,\mathbf{x})\dtau,
\end{equation}
\begin{equation}\label{if5}
 \theta(t,\mathbf{x})=  \mathcal{L}_{1}(t)\theta_{0}(\mathbf{x}) + \mathcal{L}_{2}(t)\left(\frac{1}{2}(-\Delta)\theta_{0}(\mathbf{x}) + \theta_1(\mathbf{x})\right)
+ \int_{0}^{t}\mathcal{L}_{2}(t-\tau)f_{2}(\tau,\mathbf{x})\dtau,
\end{equation}
where
\beq\label{if9}
\boldsymbol{\omega}_{0}=(\omega_{10},\omega_{20},\omega_{30}),
\eeq
\begin{equation}\label{if6}
\theta_1(\mathbf{x}) = \p_{t}\theta(0,\mathbf{x}) = -\mathbf{u}_0\cdot\nabla\theta_0 + u_{30},
\end{equation}
\[
\mathbf{u}_0 =\mathbf{u}(0,\mathbf{x})=(u_{10}, u_{20},u_{30}),
\]
with
\[
u_{10}=\partial_2(-\Delta)^{-1}\omega_{30}-\partial_3(-\Delta)^{-1}\omega_{20}, u_{20}=\partial_3(-\Delta)^{-1}\omega_{10}-\partial_1(-\Delta)^{-1}\omega_{30},
\]
\begin{equation}\label{if7}
u_{30}=\partial_1(-\Delta)^{-1}\omega_{20}-\partial_2(-\Delta)^{-1}\omega_{10}).
\end{equation}

With the help (\ref{if4})--(\ref{if5}), we focus on energy and decay estimates of solutions to (\ref{Bou3})--(\ref{inibou2}) in the following subsection.

\subsection{Energy estimates and decay estimates of nonlinear equations}

For $m\in\mathbb{N}$, we define
\[
\mathcal{E}_{1}(t) = \sup_{0\leq \tau\leq t}\left\{\|\theta(\tau)\|_{H^{m+1}}+\|\boldsymbol{\omega}(\tau)\|_{H^{m}}+\|\Lambda^{-1}\omega_{3}(\tau)\|_{L^{2}}\right\},
\]
\[
\mathcal{E}_{2}(t) = \sup_{0\leq \tau\leq t}\{\langle\tau\rangle^{\frac{5}{4}}\left(\|\nabla_{h}\theta(\tau)\|_{L^{\infty}}
+\|\nabla\mathbf{u}(\tau)\|_{L^{\infty}} +\|\boldsymbol{\omega}_h(\tau)\|_{L^{\infty}}\right)
\]
\[
+\langle\tau\rangle^{\frac{3}{2}}\|\omega_{3}(\tau)\|_{L^{\infty}}+\langle\tau\rangle\left(\|\theta(\tau)\|_{L^{\infty}}+\|\partial_{3}\theta(\tau)\|_{L^{\infty}}+\|\mathbf{u}(\tau)\|_{L^{\infty}}\right)\},
\]
\begin{equation*}
\begin{split}
  \mathcal{E}_{3}(t) =& \sup_{0\leq \tau\leq t}\left\{\langle\tau\rangle^{\frac{1}{2}}\left(\|\theta(\tau)\|_{H^{5}}+\|\mathbf{u}_{h}(\tau)\|_{H^{4}}\right)
+\langle\tau\rangle\left(\|\nabla_{h}\theta(\tau)\|_{H^{3}}+\|\boldsymbol{\omega}(\tau)\|_{H^{3}}\right)\right. \\
  &\left. + \langle\tau\rangle^{\frac{5}{4}}\left(\|\nabla_{h}^{2}\theta(\tau)\|_{H^{1}}+\|\nabla_{h}\boldsymbol{\omega}_h(\tau)\|_{H^{1}}+\|u_{3}(\tau)\|_{H^{3}}\right)\right\},
\end{split}
\end{equation*}
\begin{equation*}
  \mathcal{E}_4(t) =  \sup_{0\leq \tau\leq t}\left\{\langle\tau\rangle^{\frac{5}{4}}\left(\|\partial_{\tau}\theta(\tau)\|_{H^{1}}
  + \|\partial_{\tau}\boldsymbol{\omega}(\tau)\|_{L^{2}}
+\|\partial_{\tau}\mathbf{u}(\tau)\|_{H^{1}}\right)\right\},
\end{equation*}
\[
\mathcal{E}_{0} = \|\theta_{0}\|_{W^{10,1}} + \|\theta_{0}\|_{H^{m+1}} +\|\boldsymbol{\omega}_{0}\|_{W^{7,1}} + \|\boldsymbol{\omega}_{0}\|_{H^{m}}+ \|\Lambda^{-1}\omega_{30}\|_{L^{1}}+ \|\Lambda^{-1}\omega_{30}\|_{L^{2}},
\]
\[
\mathcal{E}(t) = \mathcal{E}_{1}(t) +\mathcal{E}_{2}(t) +\mathcal{E}_{3}(t)+\mathcal{E}_{4}(t).
\]
\begin{Lemma}\label{ELD1}
Let $m\geq 2$. Then
\begin{equation}\label{Ee1}
\mathcal{E}_{1}(t) \ls \mathcal{E}_{0}+\mathcal{E}^{\frac{3}{2}}(t)+\mathcal{E}^{2}(t).
\end{equation}
\end{Lemma}
\proof
Testing (\ref{Bou3})$_{1}$ by $\boldsymbol{\omega}$ and multiplying (\ref{Bou3})$_{2}$ by $-\Delta\theta$ give
\begin{equation}\label{e1}
 \frac{1}{2}\frac{\rm{d}}{\dt}\left( \|\boldsymbol{\omega}\|_{L^{2}}^{2}+ \|\nabla\theta\|_{L^{2}}^{2} \right) + \|\nabla\boldsymbol{\omega}\|_{L^{2}}^{2}   \\
\lesssim \|\nabla\mathbf{u}\|_{L^{\infty}}(\|\nabla\theta\|_{L^{2}}^{2}+\|\boldsymbol{\omega}\|_{L^{2}}^{2}),
\end{equation}
where we use the fact that
\[
\langle \partial_2\theta,\omega_{1} \rangle + \langle -\partial_1\theta,\omega_{2} \rangle+ \langle -u_3, -\Delta\theta\rangle
\]
\[
=\langle \partial_2\theta,\omega_{1} \rangle + \langle -\partial_1\theta,\omega_{2} \rangle+\langle  -\p_{1}(-\Delta)^{-1}\omega_{2} + \p_{2}(-\Delta)^{-1}\omega_{1},-\Delta\theta\rangle= 0.
\]
Similarly,
\[
 \langle\partial_{t}\partial^{m}\boldsymbol{\omega}, \partial^{m}\boldsymbol{\omega}\rangle + \langle\partial^{m}(\mathbf{u}\cdot \nabla \boldsymbol{\omega}), \partial^{m}\boldsymbol{\omega}\rangle-\langle\partial^{m}(\boldsymbol{\omega}\cdot \nabla \mathbf{u}), \partial^{m}\boldsymbol{\omega}\rangle - \langle\Delta\partial^{m}\boldsymbol{\omega}, \partial^{m}\boldsymbol{\omega}\rangle
\]
\begin{equation}\label{E2}
= \langle\partial^{m}\partial_{2}\theta, \partial^{m}\omega_{1}\rangle + \langle-\partial^{m}\partial_{1}\theta, \partial^{m}\omega_{2}\rangle
\end{equation}
and
\begin{equation}\label{E3}
 \langle\partial_{t}\partial^{m}\nabla\theta, \partial^{m}\nabla\theta\rangle + \langle\partial^{m}\nabla(\mathbf{u}\cdot \nabla \theta), \partial^{m}\nabla\theta\rangle
= \langle-\partial^{m}\nabla u_{3}, \partial^{m}\nabla\theta\rangle.
\end{equation}
Summing up (\ref{E2}) and (\ref{E3}) yields
\[
\frac{1}{2}\frac{\rm{d}}{\dt}\left(\|\partial^{m}\boldsymbol{\omega}\|_{L^{2}}^{2} + \|\partial^{m+1}\theta\|_{L^{2}}^{2} \right) + \|\partial^{m+1}\boldsymbol{\omega}\|_{L^{2}}^{2} =
B_{1}+B_{2}+B_{3}+B_{4}
\]
with
\[
B_{1} = - \langle\partial^{m}(\mathbf{u}\cdot\nabla\boldsymbol{\omega}), \partial^{m}\boldsymbol{\omega}\rangle,
\]
\[
B_{2} = \langle\partial^{m}(\boldsymbol{\omega}\cdot\nabla\mathbf{u}), \partial^{m}\boldsymbol{\omega}\rangle,
\]
\[
B_{3} = -\langle\partial^{m}\nabla(\mathbf{u}\cdot\nabla\theta), \partial^{m}\nabla\theta\rangle,
\]
\[
B_{4} = \langle\partial^{m}\partial_{2}\theta, \partial^{m}\omega_{1}\rangle+ \langle-\partial^{m}\partial_{1}\theta,\partial^{m}\omega_{2}\rangle
+\langle-\partial^{m}\nabla u_{3}, \partial^{m}\nabla\theta\rangle.
\]
By the commutator estimate (\ref{BL14}) in Lemma \ref{PL1},
\[
 |B_{1}|\lesssim  |\langle\partial^{m}(\mathbf{u}\cdot\nabla\boldsymbol{\omega})-\mathbf{u}\cdot\nabla\partial^{m}\boldsymbol{\omega}, \partial^{m}\boldsymbol{\omega}\rangle|
 \]
 \[
 \ls \|\langle\partial^{m}\text{ div }(\mathbf{u}\otimes\boldsymbol{\omega})-\mathbf{u}\cdot\nabla\partial^{m}\boldsymbol{\omega}\|_{L^2}\|\partial^{m}\boldsymbol{\omega}\|_{L^2}
 \]
  \begin{equation}\label{Ea2}
 \ls \|\nabla\mathbf{u}\|_{L^{\infty}}\|\boldsymbol{\omega}\|_{H^{m}}^{2} + \|\boldsymbol{\omega}\|_{L^{\infty}}\|\mathbf{u}\|_{H^{m+1}}\|\boldsymbol{\omega}\|_{H^{m}},
\end{equation}
where we use the fact that
\[
\langle\mathbf{u}\cdot\nabla\partial^{m}\boldsymbol{\omega}, \partial^{m}\boldsymbol{\omega}\rangle = 0.
\]
According to (\ref{BL11}) in Lemma \ref{PL1}, one gets
\begin{equation}\label{Ea3}
 |B_{2}| \ls \|\nabla\mathbf{u}\|_{L^{\infty}}\|\boldsymbol{\omega}\|_{H^{m}}^{2} + \|\boldsymbol{\omega}\|_{L^{\infty}}\|\nabla\mathbf{u}\|_{H^{m}}\|\boldsymbol{\omega}\|_{H^{m}}.
\end{equation}
For $B_3$, we use similar arguments to obtain
\[
|B_3|
\ls |\langle\partial^{m}\nabla(\mathbf{u}\cdot\nabla\theta)-\mathbf{u}\cdot\nabla\partial^{m}\nabla\theta, \partial^{m}\nabla\theta\rangle|
\]
\beq\label{Ea4}
\ls \|\nabla\theta\|_{L^{\infty}}\|\nabla\mathbf{u}\|_{H^{m}}\|\theta\|_{H^{m+1}} + \|\nabla\mathbf{u}\|_{L^{\infty}}\|\theta\|_{H^{m+1}}^{2},
\eeq
where we use
\[
\langle\mathbf{u}\cdot\nabla\partial^{m}\nabla\theta, \partial^{m}\nabla\theta\rangle =0.
\]
Using the fact $u_{3}=\p_{1}(-\Delta)^{-1}\omega_{2} - \p_{2}(-\Delta)^{-1}\omega_{1}$ and integrating by parts,
\[
B_{4} = -\langle\partial^{m}\nabla\p_{1}(-\Delta)^{-1}\omega_{2} , \partial^{m}\nabla\theta\rangle
+ \langle\partial^{m}\nabla\p_{2}(-\Delta)^{-1}\omega_{1}, \partial^{m}\nabla\theta\rangle
\]
\beq\label{Ea11}
+ \langle\partial^{m}\partial_{2}\theta, \partial^{m}\omega_{1}\rangle + \langle-\partial^{m}\partial_{1}\theta,\partial^{m}\omega_{2}\rangle = 0.
\eeq
Summing up (\ref{Ea2})--(\ref{Ea11}), then we apply Young's inequality and Lemma \ref{corollary1} to obtain
\[
\frac{1}{2}\frac{\rm{d}}{\dt}\left(\|\partial^{m+1}\theta\|_{L^{2}}^{2} + \|\partial^{m}\boldsymbol{\omega}\|_{L^{2}}^{2}\right) + \|\partial^{m+1}\boldsymbol{\omega}\|_{L^{2}}^{2}
\leq  C \|\boldsymbol{\omega}\|_{L^{\infty}}\|\mathbf{u}\|_{H^{m+1}}\|\boldsymbol{\omega}\|_{H^{m}}
\]
\[
+ C\left(\|\nabla\mathbf{u}\|_{L^{\infty}}+\|\boldsymbol{\omega}\|_{L^{\infty}}\right)(\|\boldsymbol{\omega}\|_{H^{m}}^{2}+\|\theta\|_{H^{m}}^{2})
+ C\|\nabla\theta\|_{L^{\infty}}\|\nabla\mathbf{u}\|_{H^{m}}\|\theta\|_{H^{m+1}}
\]
\[
\leq C \left(\|\nabla\mathbf{u}\|_{L^{\infty}} + \|\boldsymbol{\omega}\|_{L^{\infty}}\right)(\|\Lambda^{-1}\omega_{3}\|_{L^{2}}^{2}+\|\boldsymbol{\omega}\|_{H^{m}}^{2}+\|\theta\|_{H^{m+1}}^{2})
\]
\beq\label{Ea5}
+ 2C^{2}\|\nabla\theta\|_{L^{\infty}}^{2}\|\theta\|_{H^{m+1}}^{2} + \frac{1}{2}\|\nabla\mathbf{u}\|_{H^{m}}^{2},
\eeq
where $C$ is a constant. Note that
\[
\partial_{t}\omega_{3} -\Delta\omega_{3} +\mathbf{u}\cdot \nabla \omega_{3} - \boldsymbol{\omega}\cdot \nabla u_{3} = 0.
\]
Multiplying this equation by a test function $(-\Delta)^{-1}\omega_{3}$ and integrating in space yield
\[
\langle\partial_{t}\omega_{3},(-\Delta)^{-1}\omega_{3}\rangle -\langle\Delta\omega_{3},(-\Delta)^{-1}\omega_{3}\rangle +\langle\mathbf{u}\cdot \nabla\omega_{3},(-\Delta)^{-1}\omega_{3}\rangle  - \langle\boldsymbol{\omega}\cdot \nabla u_{3},(-\Delta)^{-1}\omega_{3}\rangle = 0.
\]
By using the fact that $\varphi_{3} = (-\Delta_{N})^{-1}\omega_{3}$,
\[
\frac{1}{2}\frac{\rm{d}}{\dt}\|\nabla \varphi_{3}\|_{L^{2}}^{2} + \|\Delta \varphi_{3}\|_{L^{2}}^{2} + \langle\mathbf{u}\cdot \nabla (-\Delta \varphi_{3}), \varphi_{3}\rangle - \langle\boldsymbol{\omega}\cdot \nabla u_{3}, \varphi_{3}\rangle = 0.
\]
Integrating by parts and applying Lemma \ref{PoL2} give
\[
\langle\mathbf{u}\cdot \nabla (-\Delta \varphi_{3}), \varphi_{3}\rangle - \langle\boldsymbol{\omega}\cdot \nabla u_{3}, \varphi_{3}\rangle
= \langle\mathbf{u}\Delta \varphi_{3}, \nabla \varphi_{3}\rangle + \langle\boldsymbol{\omega} u_{3}, \nabla \varphi_{3}\rangle
\]
\[
\ls \|\mathbf{u}\|_{L^{2}}\|\Delta \varphi_{3}\|_{L^{\infty}}\|\nabla \varphi_{3}\|_{L^{2}} + \|\boldsymbol{\omega}\|_{L^{\infty}}\|\nabla \varphi_{3}\|_{L^{2}}\|u_{3}\|_{L^{2}}
\]
\[
\ls \|\mathbf{u}\|_{L^{2}}\|\boldsymbol{\omega}\|_{L^{\infty}}\|\Lambda^{-1}\omega_{3}\|_{L^{2}}.
\]
Therefore,
\begin{equation}\label{Ea8}
\frac{1}{2}\frac{\rm{d}}{\dt}\|\Lambda^{-1}\omega_{3}\|_{L^{2}}^{2} + \|\omega_{3}\|_{L^{2}}^{2} \ls
 \|\mathbf{u}\|_{L^{2}}\|\boldsymbol{\omega}\|_{L^{\infty}}\|\Lambda^{-1}\omega_{3}\|_{L^{2}}.
\end{equation}
From (\ref{e1}) and (\ref{Ea5})--(\ref{Ea8}) it follows that
\[
\frac{1}{2}\frac{\rm{d}}{\dt}\left(\|\theta\|_{H^{m+1}}^{2} + \|\boldsymbol{\omega}\|_{H^{m}}^{2}+\|\Lambda^{-1}\omega_{3}\|_{L^{2}}^{2}\right) + \frac{1}{2}\|\boldsymbol{\omega}\|_{H^{m+1}}^{2}
\]
\[
\ls \left(\|\nabla\mathbf{u}\|_{L^{\infty}} + \|\boldsymbol{\omega}\|_{L^{\infty}}\right)\|\mathbf{u}\|_{H^{m+1}}\|\boldsymbol{\omega}\|_{H^{m}}
+ \|\nabla\mathbf{u}\|_{L^{\infty}}(\|\boldsymbol{\omega}\|_{H^{m}}^{2}+\|\theta\|_{H^{m+1}}^{2})
\]
\[
+ \|\nabla\theta\|_{L^{\infty}}^{2}\|\theta\|_{H^{m+1}}^{2}
+ \|\boldsymbol{\omega}\|_{L^{\infty}}\|\mathbf{u}\|_{L^{2}}\|\Lambda^{-1}\omega_{3}\|_{L^{2}}
\]
\[
\ls \left(\|\nabla\mathbf{u}\|_{L^{\infty}}+\|\boldsymbol{\omega}\|_{L^{\infty}}+\|\nabla\theta\|_{L^{\infty}}^{2})(\|\boldsymbol{\omega}\|_{H^{m}}^{2}+\|\theta\|_{H^{m+1}}^{2} +\|\Lambda^{-1}\omega_{3}\|_{L^{2}}^{2}\right).
\]
By integrating in time from 0 to $t$ we find
\[
\|\theta(t)\|_{H^{m+1}}^{2} + \|\boldsymbol{\omega}(t)\|_{H^{m}}^{2}+\|\Lambda^{-1}\omega_{3}(t)\|_{L^{2}}^{2} + \int_{0}^{t}\|\boldsymbol{\omega}(\tau)\|_{H^{m+1}}^{2}\dtau
\]
\[
\ls \|\theta_{0}\|_{H^{m+1}}^{2} + \|\boldsymbol{\omega}_{0}\|_{H^{m}}^{2}+\|\Lambda^{-1}\omega_{30}\|_{L^{2}}^{2}
 \]
\[
+ \int_{0}^{t}\langle\tau\rangle^{-\frac{5}{4}}\langle\tau\rangle^{\frac{5}{4}}\left(\|\nabla\mathbf{u}(\tau)\|_{L^{\infty}}+\|\boldsymbol{\omega}(\tau)\|_{L^{\infty}}\right)\dtau\mathcal{E}_{1}^{2}(t)
+ \int_{0}^{t}\langle\tau\rangle^{-2}\langle\tau\rangle^{2}\|\nabla\theta(\tau)\|_{L^{\infty}}^{2}\dtau\mathcal{E}_{1}^{2}(t)
\]
\[
\ls \mathcal{E}_{0}^{2} + \int_{0}^{t}\langle\tau\rangle^{-\frac{5}{4}}\dtau\mathcal{E}_{2}(t)\mathcal{E}_{1}^{2}(t)+\int_{0}^{t}\langle\tau\rangle^{-2}\dtau\mathcal{E}_{2}^{2}(t)\mathcal{E}_{1}^{2}(t)
\ls \mathcal{E}_{0}^{2} + \mathcal{E}^{3}(t)+\mathcal{E}^{4}(t).
\]
This completes the proof of Lemma \ref{ELD1}.
\hfill$\square$

\begin{Lemma}\label{ELD2}
Let $m\geq 31$. Then
\begin{equation}\label{Ea10}
\mathcal{E}_{2}(t)\ls \mathcal{E}_{0} + \mathcal{E}_{0}^{2}+\mathcal{E}^{2}(t)+\mathcal{E}^{\frac{12}{5}}(t).
\end{equation}
\end{Lemma}
\proof
From (\ref{if5})--(\ref{if7}) and Lemma \ref{AL1}--\ref{AL2}, one gets
\[
\|\nabla_{h}\theta(t)\|_{L^{\infty}} \lesssim \|\widehat{\nabla_{h}\theta}(t)\|_{\widehat{L}^{1}}
\]
\[
\ls
\langle t\rangle^{-\frac{3}{2}}\|(\xi^2+\eta^2+\pi^2k^2)^{4}\widehat{\theta_{1}}(t)\|_{\widehat{L}^{\infty}} + \langle t\rangle^{-\frac{3}{2}}\|\theta_{0}\|_{W^{10,1}}
+ I_{1}+ \cdots +I_{7}
\]
\[
\ls\langle t\rangle^{-\frac{3}{2}}\|\boldsymbol{\omega}_{0}\|_{W^{7,1}} + \langle t\rangle^{-\frac{3}{2}}\|\theta_{0}\|_{W^{10,1}} +
\langle t\rangle^{-\frac{3}{2}}\|\theta_{0}\|_{H^{m}}
(\|\boldsymbol{\omega}_{0}\|_{H^{m}}+\|\Lambda^{-1}\omega_{30}\|_{L^{2}})
\]
\begin{equation}\label{Du1}
+ I_{1}+ I_{2}+ I_{3}+ I_{4}+ I_{5}+ I_{6}+ I_{7},
\end{equation}
where
\[
I_{1} = \int_{0}^{t}\langle t-\tau \rangle^{-\frac{3}{2}}\|\partial_{\tau}\mathbf{u}(\tau)\cdot\nabla\theta(\tau)\|_{W^{8,1}}d\tau,\,
I_{2} =  \int_{0}^{t}\langle t-\tau\rangle^{-\frac{3}{2}}\|\mathbf{u}(\tau)\cdot\nabla\partial_{\tau}\theta(\tau)\|_{W^{8,1}}\dtau,
\]
\[
I_{3} =  \int_{0}^{t}\langle t-\tau\rangle^{-\frac{3}{2}}\|\mathbf{u}(\tau)\cdot\nabla\theta(\tau)\|_{W^{10,1}}\dtau,\,
I_{4} = \int_{0}^{t}\langle t-\tau\rangle^{-\frac{3}{2}}\|\partial_{1}(\mathbf{u}\cdot\nabla\omega_{2})(\tau)\|_{W^{6,1}}\dtau,
\]
\[
I_{5} =
\int_{0}^{t}\langle t-\tau\rangle^{-\frac{3}{2}}\|\partial_{2}(\mathbf{u}\cdot\nabla\omega_{1})(\tau)\|_{W^{6,1}}\dtau,\,
I_{6} =  \int_{0}^{t}\langle t-\tau\rangle^{-\frac{3}{2}}\|\partial_{1}(\boldsymbol{\omega}\cdot\nabla u_{2})(\tau)\|_{W^{6,1}}\dtau,
\]
\[
I_{7} = \int_{0}^{t}\langle t-\tau\rangle^{-\frac{3}{2}}\|\partial_{2}(\boldsymbol{\omega}\cdot\nabla u_{1})(\tau)\|_{W^{6,1}}\dtau.
\]

According to the interpolation inequality (\ref{BLT1}) together with (\ref{BL13}) in Lemma \ref{PL1m}, we have
\[
\|\partial_{t}\mathbf{u}\cdot\nabla\theta\|_{W^{8,1}}
\ls \|\partial_{t}\mathbf{u}\|_{L^{2}}\|\nabla\theta\|_{H^{8}} + \|\partial_{t}\mathbf{u}\|_{H^{8}}\|\nabla\theta\|_{L^{2}}
\]
\[
\ls \|\partial_{t}\mathbf{u}\|_{L^{2}}\|\nabla\theta\|_{H^{16}}^{\frac{1}{2}}\|\nabla\theta\|_{L^{2}}^{\frac{1}{2}} + \|\partial_{t}\mathbf{u}\|_{H^{20}}^{\frac{2}{5}}\|\partial_{t}\mathbf{u}\|_{L^{2}}^{\frac{3}{5}}\|\nabla\theta\|_{L^{2}}.
\]
From (\ref{Bou3})$_{1}$, it follows that
\[
\|\partial_{t}\mathbf{u}\|_{H^{20}}\ls \|\partial_{t}\boldsymbol{\omega}\|_{H^{19}} + \|\partial_{t}\Lambda^{-1}\omega_{3}\|_{L^{2}}
\]
\[
\ls \|\boldsymbol{\omega}\|_{H^{21}} + \|\theta\|_{H^{20}} + \|\boldsymbol{u}\|_{H^{20}}\|\boldsymbol{\omega}\|_{H^{20}}
\]
\[
\ls \|\boldsymbol{\omega}\|_{H^{21}} + \|\theta\|_{H^{20}} + \|\boldsymbol{\omega}\|_{H^{21}}^2+\|\Lambda^{-1}\omega_{3}\|_{L^{2}}\|\boldsymbol{\omega}\|_{H^{21}}.
\]
Thus
\[
I_{1}\ls \int_{0}^{t}\langle t-\tau\rangle^{-\frac{3}{2}}\langle \tau\rangle^{-\frac{3}{2}}\dtau\mathcal{E}_{4}(t)\mathcal{E}_{1}^{\frac{1}{2}}(t)\mathcal{E}_{3}^{\frac{1}{2}}(t)
\]
\[
+ \int_{0}^{t}\langle t-\tau\rangle^{-\frac{3}{2}}
\langle \tau\rangle^{-\frac{5}{4}}\dtau\mathcal{E}_{1}^{\frac{2}{5}}(t)\mathcal{E}_{4}^{\frac{3}{5}}(t)\mathcal{E}_{3}(t)
+ \int_{0}^{t}\langle t-\tau\rangle^{-\frac{3}{2}}
\langle \tau\rangle^{-\frac{5}{4}}\dtau\mathcal{E}_{1}^{\frac{4}{5}}(t)\mathcal{E}_{4}^{\frac{3}{5}}(t)\mathcal{E}_{3}(t)
\]
\beq\label{Ea15}
\ls \langle t\rangle^{-\frac{5}{4}}\left(\mathcal{E}^{2}(t) + \mathcal{E}^{\frac{12}{5}}(t)\right),
\eeq
where we use (\ref{BL22}) in Lemma \ref{PL3}.

From (\ref{BLT1}) and (\ref{BL13}) in Lemma \ref{PL1m}, one infers that
\[
\|\mathbf{u}\cdot\nabla\partial_{t}\theta\|_{W^{8,1}}\ls
\|\mathbf{u}\|_{L^{2}}\|\partial_{t}\nabla\theta\|_{H^{8}} + \|\mathbf{u}\|_{H^{8}}\|\partial_{t}\nabla\theta\|_{L^{2}}
\]
\[
\ls \|\mathbf{u}\|_{L^{2}}\|\partial_{t}\nabla\theta\|_{H^{20}}^{\frac{2}{5}}\|\partial_{t}\nabla\theta\|_{L^{2}}^{\frac{3}{5}} + \|\mathbf{u}\|_{H^{16}}^{\frac{1}{2}}\|\mathbf{u}\|_{L^{2}}^{\frac{1}{2}}\|\partial_{t}\nabla\theta\|_{L^{2}}.
\]
By applying (\ref{Bou3})$_{2}$ and Lemma \ref{corollary1},
\[
\|\partial_{t}\nabla\theta\|_{H^{20}}\ls \|u_3\|_{H^{21}} + \|\mathbf{u}\cdot\nabla\theta\|_{H^{21}}
\]
\[
\ls \|\boldsymbol{\omega}\|_{H^{20}} + \|\boldsymbol{u}\|_{H^{21}}\|\theta\|_{H^{22}}\ls \|\boldsymbol{\omega}\|_{H^{20}} + \|\boldsymbol{\omega}\|_{H^{20}}\|\theta\|_{H^{22}}+\|\Lambda^{-1}\omega_{3}\|_{L^{2}}\|\theta\|_{H^{22}}.
\]
Then
\[
I_{2}\ls \int_{0}^{t}\langle t-\tau\rangle^{-\frac{3}{2}}
\langle\tau\rangle^{-\frac{5}{4}}\dtau\mathcal{E}_{3}(t)\mathcal{E}_{1}^{\frac{2}{5}}(t)\mathcal{E}_{4}^{\frac{3}{5}}(t)
\]
\[
+ \int_{0}^{t}\langle t-\tau\rangle^{-\frac{3}{2}}
\langle\tau\rangle^{-\frac{5}{4}}\dtau\mathcal{E}_{3}(t)\mathcal{E}_{1}^{\frac{4}{5}}(t)\mathcal{E}_{4}^{\frac{3}{5}}(t)
+ \int_{0}^{t}\langle t-\tau\rangle^{-\frac{3}{2}}
\langle \tau\rangle^{-\frac{3}{2}}\dtau\mathcal{E}_{3}^{\frac{1}{2}}(t)\mathcal{E}_{1}^{\frac{1}{2}}(t)\mathcal{E}_{4}(t)
\]
\beq\label{Ea16}
\ls \langle t\rangle^{-\frac{5}{4}}\left(\mathcal{E}^{2}(t) + \mathcal{E}^{\frac{12}{5}}(t)\right).
\eeq

By the interpolation inequality (\ref{BLT1}) and (\ref{BL13}) in Lemma \ref{PL1m},
\[
\|\mathbf{u}\cdot\nabla\theta\|_{W^{10,1}}
\ls \|\mathbf{u}_{h}\|_{L^{2}}\|\nabla_{h}\theta\|_{H^{10}} + \|\mathbf{u}_{h}\|_{H^{10}}\|\nabla_{h}\theta\|_{L^{2}}
\]
\[
+ \|u_{3}\|_{L^{2}}\|\partial_{3}\theta\|_{H^{10}} + \|u_{3}\|_{H^{10}}\|\partial_{3}\theta\|_{L^{2}}
\]
\[
\ls \|\mathbf{u}_{h}\|_{L^{2}}\|\nabla_{h}\theta\|_{H^{31}}^{\frac{1}{4}}\|\nabla_{h}\theta\|_{H^{3}}^{\frac{3}{4}} + \|\mathbf{u}_{h}\|_{H^{20}}^{\frac{1}{2}}\|\mathbf{u}_{h}\|_{L^{2}}^{\frac{1}{2}}\|\nabla_{h}\theta\|_{L^{2}}
\]
\[
+ \|u_{3}\|_{L^{2}}\|\partial_{3}\theta\|_{H^{20}}^{\frac{1}{2}}\|\partial_{3}\theta\|_{L^{2}}^{\frac{1}{2}} + \|u_{3}\|_{H^{25}}^{\frac{2}{5}}\|u_{3}\|_{L^{2}}^{\frac{3}{5}}\|\partial_{3}\theta\|_{L^{2}}.
\]
This implies that for $m\geq31$,
\[
I_{3}\ls \int_{0}^{t}\langle t-\tau\rangle^{-\frac{3}{2}} \langle\tau\rangle^{-\frac{5}{4}}\dtau\mathcal{E}_{3}^{\frac{7}{4}}(t)\mathcal{E}_{1}^{\frac{1}{4}}(t)
\]
\[
+ \int_{0}^{t}\langle t-\tau\rangle^{-\frac{3}{2}}\langle \tau\rangle^{-\frac{5}{4}}\dtau\mathcal{E}_{3}^{\frac{3}{2}}(t)\mathcal{E}_{1}^{\frac{1}{2}}(t)
+\int_{0}^{t}\langle t-\tau\rangle^{-\frac{3}{2}}\langle \tau\rangle^{-\frac{3}{2}}\dtau\mathcal{E}_{3}^{\frac{3}{2}}(t)\mathcal{E}_{1}^{\frac{1}{2}}(t)
\]
\beq\label{Ea18}
+\int_{0}^{t}\langle t-\tau\rangle^{-\frac{3}{2}}
\langle \tau\rangle^{-\frac{5}{4}}\dtau\mathcal{E}_{3}^{\frac{8}{5}}(t)\mathcal{E}_{1}^{\frac{2}{5}}(t)
\ls \langle t\rangle^{-\frac{5}{4}}\mathcal{E}^{2}(t).
\eeq

For $I_{4}$, we have
\[
\|\partial_{1}(\mathbf{u}\cdot\nabla\omega_{2})\|_{W^{6,1}}
\ls \|\mathbf{u}\cdot\nabla\omega_{2}\|_{W^{7,1}}
\ls \|\mathbf{u}\|_{H^{8}}\|\omega_{2}\|_{L^{2}}+\|\mathbf{u}\|_{L^{2}}\|\omega_{2}\|_{H^{8}}
\]
\[
\ls \|\mathbf{u}\|_{L^{2}}^{\frac{2}{3}}\|\mathbf{u}\|_{H^{24}}^{\frac{1}{3}}\|\omega_{2}\|_{L^{2}}+\|\mathbf{u}\|_{L^{2}}\|\omega_{2}\|_{H^{23}}^{\frac{1}{4}}\|\omega_{2}\|_{H^{3}}^{\frac{3}{4}}.
\]
Then
\[
I_{4}\ls \int_{0}^{t}\langle t-\tau\rangle^{-\frac{3}{2}}\langle \tau\rangle^{-\frac{4}{3}}\dtau\mathcal{E}_{3}^{\frac{5}{3}}(t)\mathcal{E}_{1}^{\frac{1}{3}}(t)
\]
\beq\label{Ea21}
+ \int_{0}^{t}\langle t-\tau\rangle^{-\frac{3}{2}}\langle \tau\rangle^{-\frac{5}{4}}\dtau\mathcal{E}_{3}^{\frac{7}{4}}(t)\mathcal{E}_{1}^{\frac{1}{4}}(t)
\ls \langle t\rangle^{-\frac{5}{4}}\mathcal{E}^{2}(t).
\eeq

Similarly for $I_{5}$,
\beq\label{Ea22}
I_{5}\ls \langle t\rangle^{-\frac{5}{4}}\mathcal{E}^{2}(t).
\eeq

For $I_{6}$, we have
\[
\|\partial_{1}(\boldsymbol{\omega}\cdot\nabla u_{2})\|_{W^{6,1}}\ls \|\boldsymbol{\omega}\cdot\nabla u_{2}\|_{W^{7,1}}
\ls \|\boldsymbol{\omega}\|_{H^{8}}\|u_{2}\|_{L^{2}}+\|\boldsymbol{\omega}\|_{L^{2}}\|u_{2}\|_{H^{8}}
\]
\[
\ls \|\boldsymbol{\omega}\|_{H^{3}}^{\frac{3}{4}}\|\boldsymbol{\omega}\|_{H^{23}}^{\frac{1}{4}}\|u_{2}\|_{L^{2}}+\|\boldsymbol{\omega}\|_{L^{2}}\|u_{2}\|_{H^{16}}^{\frac{1}{2}}\|u_{2}\|_{L^{2}}^{\frac{1}{2}}.
\]
Then
\[
I_{6}\ls \int_{0}^{t}\langle t-\tau\rangle^{-\frac{3}{2}}\langle \tau\rangle^{-\frac{5}{4}}\dtau\mathcal{E}_{3}^{\frac{7}{4}}(t)\mathcal{E}_{1}^{\frac{1}{4}}(t)
\]
\beq\label{Ea24}
+ \int_{0}^{t}\langle t-\tau\rangle^{-\frac{3}{2}}\langle \tau\rangle^{-\frac{5}{4}}\dtau\mathcal{E}_{3}^{\frac{3}{2}}(t)\mathcal{E}_{1}^{\frac{1}{2}}(t)
\ls \langle t\rangle^{-\frac{5}{4}}\mathcal{E}^{2}(t).
\eeq

Similarly for $I_{7}$,
\beq\label{Ea25}
I_{7}\ls \langle t\rangle^{-\frac{5}{4}}\mathcal{E}^{2}(t).
\eeq

Substituting (\ref{Ea15})--(\ref{Ea25}) into (\ref{Du1}) gives
\beq\label{Ef1}
\|\nabla_{h}\theta(t)\|_{L^{\infty}} \lesssim \|\widehat{\nabla_{h}\theta}(t)\|_{\widehat{L}^{1}}
\ls \langle t\rangle^{-\frac{5}{4}}\left(\mathcal{E}_{0}+\mathcal{E}_{0}^{2}+\mathcal{E}^{2}(t)+\mathcal{E}^{\frac{12}{5}}(t)\right).
\eeq
Note that
\[
\|\theta(t)\|_{L^{\infty}}+\|\partial_{3}\theta(t)\|_{L^{\infty}}\lesssim
\|\widehat{\theta}(t)\|_{\widehat{L}^{1}}+\|\widehat{\partial_{3}\theta}(t)\|_{\widehat{L}^{1}}
\]
\[
\ls\langle t\rangle^{-1}\|\boldsymbol{\omega}_{0}\|_{W^{6,1}} + \langle t\rangle^{-1}\|\theta_{0}\|_{W^{9,1}} +
\langle t\rangle^{-1}\|\theta_{0}\|_{H^{m}}
(\|\boldsymbol{\omega}_{0}\|_{H^{m}}+\|\Lambda^{-1}\omega_{30}\|_{L^{2}})
\]
\[
+ \int_{0}^{t}\langle t-\tau \rangle^{-1}\|\partial_{\tau}\mathbf{u}(\tau)\cdot\nabla\theta(\tau)\|_{W^{7,1}}\dtau +
\int_{0}^{t}\langle t-\tau\rangle^{-1}\|\mathbf{u}(\tau)\cdot\nabla\partial_{\tau}\theta(\tau)\|_{W^{7,1}}\dtau
\]
\[
+\int_{0}^{t}\langle t-\tau\rangle^{-1}\|\mathbf{u}(\tau)\cdot\nabla\theta(\tau)\|_{W^{9,1}}\dtau +
\int_{0}^{t}\langle t-\tau\rangle^{-1}\|\partial_{1}(\mathbf{u}\cdot\nabla\omega_{2})(\tau)\|_{W^{5,1}}\dtau
\]
\[
+\int_{0}^{t}\langle t-\tau\rangle^{-1}\|\partial_{2}(\mathbf{u}\cdot\nabla\omega_{1})(\tau)\|_{W^{5,1}}\dtau + \int_{0}^{t}\langle t-\tau\rangle^{-1}\|\partial_{1}(\boldsymbol{\omega}\cdot\nabla u_{2})(\tau)\|_{W^{5,1}}\dtau
\]
\[
+\int_{0}^{t}\langle t-\tau\rangle^{-1}\|\partial_{2}(\boldsymbol{\omega}\cdot\nabla u_{1})(\tau)\|_{W^{5,1}}\dtau.
\]
In the similar way,
\beq\label{Ef2}
\|\theta(t)\|_{L^{\infty}}+\|\partial_{3}\theta(t)\|_{L^{\infty}}
\ls \langle t\rangle^{-1}\left(\mathcal{E}_{0}+\mathcal{E}_{0}^{2}+\mathcal{E}^{2}(t)++\mathcal{E}^{\frac{12}{5}}(t)\right).
\eeq

Now we turn to $\boldsymbol{\omega}$. From (\ref{if4}), one deduces that
\[
\|\omega_{1}(t)\|_{L^{\infty}}\ls\|\widehat{\omega_{1}}(t)\|_{\widehat{L}^{1}}
\]
\[
\lesssim \|e^{-\left(\xi^{2}+\eta^2+ \pi^2k ^{2}\right) t}\widehat{\omega_{10}}\|_{\widehat{L}^{1}}
+ \int_{0}^{t}\|e^{-\left(\xi^{2}+\eta^2+ \pi^2k ^{2}\right)(t-\tau)}\widehat{\mathbf{u}\cdot\nabla\omega_{1}}(\tau)\|_{\widehat{L}^{1}} \dtau
\]
\[
 +
\int_{0}^{t}\|e^{-\left(\xi^{2}+\eta^2+\pi^2 k ^{2}\right)(t-\tau)}\widehat{\boldsymbol{\omega}\cdot\nabla u_{1}}(\tau)\|_{\widehat{L}^{1}} \dtau
+\int_{0}^{t}\|e^{-\left(\xi^{2}+\eta^2+ \pi^2k ^{2}\right)(t-\tau)}\widehat{\p_{2}\theta}(\tau)\|_{\widehat{L}^{1}}\dtau
\]
\[
\lesssim e^{- t}\|\widehat{\boldsymbol{\omega}}_{0}\|_{\widehat{L}^{1}} +
\int_{0}^{t}e^{-(t-\tau)}\|\widehat{\mathbf{u}\cdot\nabla\omega_{1}}(\tau)\|_{\widehat{L}^{1}}\dtau
\]
\[
+ \int_{0}^{t}e^{-(t-\tau)}\|\widehat{\boldsymbol{\omega}\cdot\nabla u_{1}}(\tau)\|_{\widehat{L}^{1}}\dtau+
\int_{0}^{t}e^{-(t-\tau)}\langle\tau\rangle^{-\frac{5}{4}}\langle\tau\rangle^{\frac{5}{4}}\|\widehat{\nabla_{h}\theta}(\tau)\|_{\widehat{L}^{1}}\dtau
\]
\[
\lesssim e^{- t}\|\boldsymbol{\omega}_{0}\|_{H^{2}} +
\int_{0}^{t}e^{-(t-\tau)}\|\mathbf{u}(\tau)\cdot\nabla \omega_{1}(\tau)\|_{H^{2}}\dtau + \int_{0}^{t}e^{-(t-\tau)}\|\boldsymbol{\omega}(\tau)\cdot\nabla u_{1}(\tau)\|_{H^{2}}\dtau
\]
\beq\label{Ef3}
+ \langle t\rangle^{-\frac{5}{4}}
\left(\mathcal{E}_{0} + \mathcal{E}_{0}^{2}+\mathcal{E}^{2}(t)+\mathcal{E}^{\frac{12}{5}}(t)\right),
\eeq
where we use (\ref{eLL3}), (\ref{Ef1}) and the fact $k\geq1$. Indeed, we have
\[
\|\mathbf{u}\cdot\nabla \omega_{1}\|_{H^{2}}
\ls \|\mathbf{u}\|_{H^{3}}\|\omega_{1}\|_{H^{3}}.
\]
Then
\[
\int_{0}^{t}e^{-(t-\tau)}\|\mathbf{u}(\tau)\cdot\nabla \omega_{1}(\tau)\|_{H^{2}}\dtau
\ls \int_{0}^{t}e^{-(t-\tau)}\langle \tau\rangle^{-\frac{3}{2}}\dtau\mathcal{E}_{3}^{2}(t)
\]
\beq\label{Ef4}
\ls \langle t\rangle^{-\frac{3}{2}}\mathcal{E}^{2}(t).
\eeq
Similarly,
\beq\label{Ef5}
\int_{0}^{t}e^{-(t-\tau)}\|\boldsymbol{\omega}(\tau)\cdot\nabla u_{1}(\tau)\|_{H^{2}}\dtau
\ls \langle t\rangle^{-\frac{3}{2}}\mathcal{E}^{2}(t).
\eeq
By plugging (\ref{Ef4}) and (\ref{Ef5}) into (\ref{Ef3}),
\beq\label{Ef6}
\|\omega_{1}(t)\|_{L^{\infty}}\ls\|\widehat{\omega_{1}}(t)\|_{\widehat{L}^{1}}
\ls \langle t\rangle^{-\frac{5}{4}}
\left(\mathcal{E}_{0} + \mathcal{E}_{0}^{2}+\mathcal{E}^{2}(t)+\mathcal{E}^{\frac{12}{5}}(t)\right).
\eeq
Using bounds similar to those established for $\|\omega_{1}(t)\|_{L^{\infty}}$ gives
\beq\label{Ef7}
\|\omega_{2}(t)\|_{L^{\infty}}\ls\|\widehat{\omega_{2}}(t)\|_{\widehat{L}^{1}}
\ls \langle t\rangle^{-\frac{5}{4}}
\left(\mathcal{E}_{0} + \mathcal{E}_{0}^{2}+\mathcal{E}^{2}(t)+\mathcal{E}^{\frac{12}{5}}(t)\right).
\eeq
According to (\ref{if4}), we see that
\begin{equation}\label{Ef8}
   \widehat{\omega_{3}}(t,\xi,\eta,k) = e^{-\left(\xi^{2}+\eta^2+ \pi^2k ^{2}\right)t}\widehat{\omega_{30}}(\xi,\eta,k) + \int_{0}^{t}e^{-\left(\xi^{2}+\eta^2+ \pi^2k ^{2}\right)(t-\tau)}\widehat{f_{13}}(\tau,\xi,\eta,k)\dtau,
\end{equation}
where
$(\xi,\eta)\in\mathbb{R}\times\mathbb{R}$ and $k\geq0$. Then
\[
\|\widehat{\Lambda^{-1}\omega_{3}}(t)\|_{\widehat{L}^{1}}
\ls \|e^{-\left(\xi^{2}+\eta^2+\pi^2k ^{2}\right) t}\left(\xi^{2}+\eta^2+ \pi^{2}k ^{2}\right)^{-\frac{1}{2}}\widehat{\omega_{30}}\|_{\widehat{L}^{1}}
\]
\[
+
\int_{0}^{t}\|e^{-\left(\xi^{2}+\eta^2+ \pi^{2}k ^{2}\right)(t-\tau)}\left(\xi^{2}+\eta^2+ \pi^{2}k ^{2}\right)^{-\frac{1}{2}}\widehat{\boldsymbol{\omega}\cdot\nabla u_{3}}(\tau)\|_{\widehat{L}^{1}}\dtau
\]
\[
+\int_{0}^{t}\|e^{-\left(\xi^{2}+\eta^2+ \pi^{2}k ^{2}\right)(t-\tau)}\left(\xi^{2}+\eta^2+ \pi^{2}k ^{2}\right)^{-\frac{1}{2}}\widehat{\mathbf{u}\cdot\nabla\omega_{3}}(\tau)\|_{\widehat{L}^{1}}\dtau
\]
\beq\label{EW1}
:= U_{1}+U_{2}+U_{3}.
\eeq
Since $k\geq 0$, the first term $U_{1}$ is bounded as
\[
U_{1}
\ls \|\left(\xi^{2}+\eta^2+ \pi^{2}k ^{2}\right)^{-\frac{1}{2}}\widehat{\omega_{30}}\|_{\widehat{L}^{\infty}}\left(\int_{\mathbb{R}^{2}}e^{-\left(\xi^{2}+\eta^{2}\right)t}\dxi\deta
+ \int_{\mathbb{R}^{3}}e^{-\left(\xi^{2}+\eta^{2}+\zeta^{2}\right)t}\dxi\deta\dzeta\right)
\]
\[
\ls \left(t^{-1}+t^{-\frac{3}{2}}\right)\|\Lambda^{-1}\omega_{30}\|_{L^{1}}\ls t^{-1}\|\Lambda^{-1}\omega_{30}\|_{L^{1}}.
\]
From (\ref{eLL3}), one also deduces that
\[
U_{1}\ls \|\left(\xi^{2}+\eta^2+ \pi^{2}k ^{2}\right)^{-\frac{1}{2}}\widehat{\omega_{30}}\|_{\widehat{L}^{1}}\ls \|\Lambda^{-1}\omega_{30}\|_{H^{2}}.
\]
Therefore,
\beq\label{EK3}
U_{1}\ls \langle t\rangle^{-1} \mathcal{E}_{0}.
\eeq
Using the fact $\boldsymbol{\omega}\cdot\nabla u_{3}=\text{div}(\boldsymbol{\omega}u_{3})$,
\[
U_2\ls \int_{0}^{t}\|e^{-\left(\xi^{2}+\eta^2+ \pi^{2}k ^{2}\right)(t-\tau)}\widehat{\boldsymbol{\omega}u_{3}}(\tau)\|_{\widehat{L}^{1}}\dtau
\]
\[
\ls\int_{0}^{t}(t-\tau)^{-1}\|\boldsymbol{\omega}(\tau)u_{3}(\tau)\|_{L^{1}}\dtau
\ls\int_{0}^{t}(t-\tau)^{-1}\|\boldsymbol{\omega}(\tau)\|_{L^{2}}\|u_{3}(\tau)\|_{L^{2}}\dtau.
\]
Moreover,
\[
U_2\ls\int_{0}^{t}\|\widehat{\boldsymbol{\omega}u_{3}}(\tau)\|_{\widehat{L}^{1}}\dtau\ls\int_{0}^{t}\|\boldsymbol{\omega}(\tau)u_{3}(\tau)\|_{H^{2}}\dtau
\ls\int_{0}^{t}\|\boldsymbol{\omega}(\tau)\|_{H^{2}}\|u_{3}(\tau)\|_{H^{2}}\dtau.
\]
Thus
\beq\label{EK6}
U_2
\ls \int_{0}^{t}\langle t-\tau\rangle^{-1}\langle \tau\rangle^{-\frac{9}{4}}\dtau\mathcal{E}_{3}^{2}(t)
\ls \langle t\rangle^{-1}\mathcal{E}^{2}(t).
\eeq
In the similar way,
\[
U_3 \ls \int_{0}^{t}\langle t-\tau\rangle^{-1}\|\mathbf{u}(\tau)\|_{H^{2}}\|\omega_{3}(\tau)\|_{H^{2}}\dtau
\]
\beq\label{EK7}
\ls \int_{0}^{t}\langle t-\tau\rangle^{-1}
\langle \tau\rangle^{-\frac{3}{2}}\dtau\mathcal{E}_{3}^{2}(t)
\ls \langle t\rangle^{-1}\mathcal{E}^{2}(t).
\eeq
By substituting (\ref{EK3})--(\ref{EK7}) into (\ref{EW1}),
\beq\label{EW2}
\|\Lambda^{-1}\omega_{3}(t)\|_{L^{\infty}}\ls\|\widehat{\Lambda^{-1}\omega_{3}}(t)\|_{\widehat{L}^{1}}\ls \langle t\rangle^{-1}\left(\mathcal{E}_{0}+\mathcal{E}^{2}(t)\right).
\eeq
Note that
\[
\|\omega_{3}(t)\|_{L^{\infty}}\ls\|\widehat{\omega_{3}}(t)\|_{\widehat{L}^{1}}
\]
\[
\ls \|e^{-\left(\xi^{2}+\eta^2+ \pi^{2}k ^{2}\right) t}\left(\xi^{2}+\eta^2+ \pi^{2}k ^{2}\right)^{\frac{1}{2}}\left(\xi^{2}+\eta^2+ \pi^{2}k ^{2}\right)^{-\frac{1}{2}}\widehat{\omega_{30}}\|_{\widehat{L}^{1}}
\]
\[
+
\int_{0}^{t}\|e^{-\left(\xi^{2}+\eta^2+ \pi^{2}k ^{2}\right)(t-\tau)}\left(\xi^{2}+\eta^2+ \pi^{2}k ^{2}\right)^{\frac{1}{2}}\left(\xi^{2}+\eta^2+ \pi^{2}k ^{2}\right)^{-\frac{1}{2}}\widehat{\boldsymbol{\omega}\cdot\nabla u_{3}}(\tau)\|_{\widehat{L}^{1}}\dtau
\]
\[
 +\int_{0}^{t}\|e^{-\left(\xi^{2}+\eta^2+ \pi^{2}k ^{2}\right)(t-\tau)}\left(\xi^{2}+\eta^2+ \pi^{2}k ^{2}\right)^{\frac{1}{2}}\left(\xi^{2}+\eta^2+ \pi^{2}k ^{2}\right)^{-\frac{1}{2}}\widehat{\mathbf{u}\cdot\nabla\omega_{3}}(\tau)\|_{\widehat{L}^{1}}\dtau.
\]
Hence
\[
\|\omega_{3}(t)\|_{L^{\infty}}\ls\|\widehat{\omega_{3}}(t)\|_{\widehat{L}^{1}}
\]
\[
\ls \langle t\rangle^{-\frac{3}{2}}\left(\|\Lambda^{-1}\omega_{30}\|_{L^{1}}+\|\Lambda^{-1}\omega_{30}\|_{H^{2}}\right)+\int_{0}^{t}\langle t-\tau\rangle^{-\frac{3}{2}}\|\boldsymbol{\omega}(\tau)u_{3}(\tau)\|_{H^{3}}\dtau
\]
\[
+\int_{0}^{t}\langle t-\tau\rangle^{-\frac{3}{2}}\|\mathbf{u}(\tau)\omega_{3}(\tau)\|_{H^{3}}\dtau
\]
\[
\ls \langle t\rangle^{-\frac{3}{2}} \mathcal{E}_{0}+\int_{0}^{t}\langle t-\tau\rangle^{-\frac{3}{2}}\langle \tau\rangle^{-\frac{9}{4}}\dtau\mathcal{E}_{3}^{2}(t)
+ \int_{0}^{t}\langle t-\tau\rangle^{-\frac{3}{2}}\langle \tau\rangle^{-\frac{3}{2}}\dtau\mathcal{E}_{3}^{2}(t)
\]
\beq\label{EW6}
\ls \langle t\rangle^{-\frac{3}{2}}\left(\mathcal{E}_{0}+\mathcal{E}^{2}(t)\right),
\eeq
where we use $e^{-\frac{1}{2}\left(\xi^{2}+\eta^2+ \pi^{2}k ^{2}\right) t}\left(\xi^{2}+\eta^2+ \pi^{2}k ^{2}\right)^{\frac{1}{2}}t^{\frac{1}{2}}\ls1$. Recalling that
\[
u_{1} = \partial_{2}(-\Delta)^{-1}\omega_{3}-\partial_{3}(-\Delta)^{-1}\omega_{2},
\]
then
\[
\|u_{1}(t)\|_{L^{\infty}}
\ls\|\mathcal{F}\left(\partial_{2}(-\Delta)^{-1}\omega_{3}\right)(t)\|_{\widehat{L}^{1}}
+\|\mathcal{F}\left(\partial_{3}(-\Delta)^{-1}\omega_{2}\right)(t)\|_{\widehat{L}^{1}}.
\]
Furthermore,
\[
\|\mathcal{F}\left(\partial_{2}(-\Delta)^{-1}\omega_{3}\right)(t)\|_{\widehat{L}^{1}}\ls
\left\|\frac{i\eta}{\xi^2+\eta^2+\pi^{2} k^2}\widehat{\omega_{3}}(t)\right\|_{\widehat{L}^{1}}
\]
\[
\ls\|\widehat{\Lambda^{-1}\omega_{3}}(t)\|_{\widehat{L}^{1}},\,(\xi,\eta)\in\mathbb{R}
\times\mathbb{R},\, k\geq 0,
\]
\[
\|\mathcal{F}\left(\partial_{3}(-\Delta)^{-1}\omega_{2}\right)(t)\|_{\widehat{L}^{1}}
\ls
\left\|\frac{k}{\xi^2+\eta^2+\pi^{2} k^2}\widehat{\omega_{2}}(t)\right\|_{\widehat{L}^{1}}
\]
\[
\ls \|\widehat{\omega_{2}}(t)\|_{\widehat{L}^{1}},\,(\xi,\eta)\in\mathbb{R}
\times\mathbb{R},\,k\geq 1.
\]
Hence we use (\ref{Ef7}) and (\ref{EW2}) to get
\[
\|u_{1}(t)\|_{L^{\infty}}\ls \|\widehat{\Lambda^{-1}\omega_{3}}(t)\|_{\widehat{L}^{1}}+\|\widehat{\omega_{2}}(t)\|_{\widehat{L}^{1}}
\ls \langle t\rangle^{-1}
\left(\mathcal{E}_{0} + \mathcal{E}_{0}^{2}+\mathcal{E}^{2}(t)+\mathcal{E}^{\frac{12}{5}}(t)\right).
\]
In the similar way,
\[
\|\nabla u_{1}(t)\|_{L^{\infty}}\ls \|\widehat{\omega_{3}}(t)\|_{\widehat{L}^{1}} + \|\widehat{\omega_{2}}(t)\|_{\widehat{L}^{1}}
\ls \langle t\rangle^{-\frac{5}{4}}
\left(\mathcal{E}_{0} + \mathcal{E}_{0}^{2}+\mathcal{E}^{2}(t)+\mathcal{E}^{\frac{12}{5}}(t)\right),
\]
\[
\| u_{2}(t)\|_{L^{\infty}}\ls \|\widehat{\Lambda^{-1}\omega_{3}}(t)\|_{\widehat{L}^{1}} + \|\widehat{\omega_{1}}(t)\|_{\widehat{L}^{1}}
\ls \langle t\rangle^{-1}
\left(\mathcal{E}_{0} + \mathcal{E}_{0}^{2}+\mathcal{E}^{2}(t)+\mathcal{E}^{\frac{12}{5}}(t)\right),
\]
\[
\|\nabla u_{2}(t)\|_{L^{\infty}}\ls \|\widehat{\omega_{3}}(t)\|_{\widehat{L}^{1}} + \|\widehat{\omega_{1}}(t)\|_{\widehat{L}^{1}}
\ls \langle t\rangle^{-\frac{5}{4}}
\left(\mathcal{E}_{0} + \mathcal{E}_{0}^{2}+\mathcal{E}^{2}(t)+\mathcal{E}^{\frac{12}{5}}(t)\right),
\]
\[
\|u_{3}(t)\|_{W^{1,\infty}}\ls \|\widehat{\omega_{1}}(t)\|_{\widehat{L}^{1}}+\|\widehat{\omega_{2}}(t)\|_{\widehat{L}^{1}}
\ls \langle t\rangle^{-\frac{5}{4}}
\left(\mathcal{E}_{0} + \mathcal{E}_{0}^{2}+\mathcal{E}^{2}(t)+\mathcal{E}^{\frac{12}{5}}(t)\right).
\]
Finally, we finish the proof of Lemma \ref{ELD2}.
\hfill$\square$
\begin{Lemma}\label{ELD3}
Let $m\geq 31$. Then
\begin{equation}\label{EW12}
\mathcal{E}_{3}(t)\ls \mathcal{E}_{0} + \mathcal{E}_{0}^{2}+\mathcal{E}^{2}(t)+\mathcal{E}^{\frac{12}{5}}(t).
\end{equation}
\end{Lemma}
\proof
From (\ref{if5})--(\ref{if7}) and Lemma \ref{AL1}--\ref{AL2}, one gets
\[
\|\nabla_{h}^{2}\theta(t)\|_{H^{1}}
\ls\langle t\rangle^{-\frac{3}{2}}\|\boldsymbol{\omega}_{0}\|_{W^{7,1}} + \langle t\rangle^{-\frac{3}{2}}\|\theta_{0}\|_{W^{10,1}}
\]
\[
+
\langle t\rangle^{-\frac{3}{2}}\|\theta_{0}\|_{H^{m}}
(\|\boldsymbol{\omega}_{0}\|_{H^{m}}+\|\Lambda^{-1}\omega_{30}\|_{L^{2}})
+ \int_{0}^{t}\langle t-\tau \rangle^{-\frac{3}{2}}\|\partial_{\tau}\mathbf{u}(\tau)\cdot\nabla\theta(\tau)\|_{W^{8,1}}\dtau
\]
\[
+  \int_{0}^{t}\langle t-\tau\rangle^{-\frac{3}{2}}\|\mathbf{u}(\tau)\cdot\nabla\partial_{\tau}\theta(\tau)\|_{W^{8,1}}\dtau
+  \int_{0}^{t}\langle t-\tau\rangle^{-\frac{3}{2}}\|\mathbf{u}(\tau)\cdot\nabla\theta(\tau)\|_{W^{10,1}}\dtau
\]
\[
+ \int_{0}^{t}\langle t-\tau\rangle^{-\frac{3}{2}}\|\partial_{1}(\mathbf{u}\cdot\nabla\omega_{2})(\tau)\|_{W^{6,1}}\dtau
+
\int_{0}^{t}\langle t-\tau\rangle^{-\frac{3}{2}}\|\partial_{2}(\mathbf{u}\cdot\nabla\omega_{1})(\tau)\|_{W^{6,1}}\dtau
\]
\[
+  \int_{0}^{t}\langle t-\tau\rangle^{-\frac{3}{2}}\|\partial_{1}(\boldsymbol{\omega}\cdot\nabla u_{2})(\tau)\|_{W^{6,1}}\dtau
+ \int_{0}^{t}\langle t-\tau\rangle^{-\frac{3}{2}}\|\partial_{2}(\boldsymbol{\omega}\cdot\nabla u_{1})(\tau)\|_{W^{6,1}}\dtau.
\]
Similar to the $L^\infty$-estimate of $\nabla_h \theta$ in Lemma \ref{ELD2},
\beq\label{EN1}
\|\nabla_{h}^{2}\theta(t)\|_{H^{1}}
\ls \langle t\rangle^{-\frac{5}{4}}\left(\mathcal{E}_{0} + \mathcal{E}_{0}^{2}+\mathcal{E}^{2}(t)+\mathcal{E}^{\frac{12}{5}}(t)\right).
\eeq
Moreover,
\[
\|\nabla_{h}\theta(t)\|_{H^{3}}
\ls\langle t\rangle^{-1}\|\boldsymbol{\omega}_{0}\|_{W^{7,1}} + \langle t\rangle^{-1}\|\theta_{0}\|_{W^{10,1}}
\]
\[
 +
\langle t\rangle^{-1}\|\theta_{0}\|_{H^{m}}
(\|\boldsymbol{\omega}_{0}\|_{H^{m}}+\|\Lambda^{-1}\omega_{30}\|_{L^{2}})
+ \int_{0}^{t}\langle t-\tau \rangle^{-1}\|\partial_{\tau}\mathbf{u}(\tau)\cdot\nabla\theta(\tau)\|_{W^{8,1}}\dtau
\]
\[
+  \int_{0}^{t}\langle t-\tau\rangle^{-1}\|\mathbf{u}(\tau)\cdot\nabla\partial_{\tau}\theta(\tau)\|_{W^{8,1}}\dtau
+  \int_{0}^{t}\langle t-\tau\rangle^{-1}\|\mathbf{u}(\tau)\cdot\nabla\theta(\tau)\|_{W^{10,1}}\dtau
\]
\[
+ \int_{0}^{t}\langle t-\tau\rangle^{-1}\|\partial_{1}(\mathbf{u}\cdot\nabla\omega_{2})(\tau)\|_{W^{6,1}}\dtau
+
\int_{0}^{t}\langle t-\tau\rangle^{-1}\|\partial_{2}(\mathbf{u}\cdot\nabla\omega_{1})(\tau)\|_{W^{6,1}}\dtau
\]
\[
+  \int_{0}^{t}\langle t-\tau\rangle^{-1}\|\partial_{1}(\boldsymbol{\omega}\cdot\nabla u_{2})(\tau)\|_{W^{6,1}}\dtau
+ \int_{0}^{t}\langle t-\tau\rangle^{-1}\|\partial_{2}(\boldsymbol{\omega}\cdot\nabla u_{1})(\tau)\|_{W^{6,1}}\dtau,
\]
\[
\|\theta(t)\|_{H^{5}}
\ls\langle t\rangle^{-\frac{1}{2}}\|\boldsymbol{\omega}_{0}\|_{W^{7,1}} + \langle t\rangle^{-\frac{1}{2}}\|\theta_{0}\|_{W^{10,1}}
\]
\[
 +
\langle t\rangle^{-\frac{1}{2}}\|\theta_{0}\|_{H^{m}}
(\|\boldsymbol{\omega}_{0}\|_{H^{m}}+\|\Lambda^{-1}\omega_{30}\|_{L^{2}})
+ \int_{0}^{t}\langle t-\tau \rangle^{-\frac{1}{2}}\|\partial_{\tau}\mathbf{u}(\tau)\cdot\nabla\theta(\tau)\|_{W^{8,1}}\dtau
\]
\[
+  \int_{0}^{t}\langle t-\tau\rangle^{-\frac{1}{2}}\|\mathbf{u}(\tau)\cdot\nabla\partial_{\tau}\theta(\tau)\|_{W^{8,1}}\dtau
+  \int_{0}^{t}\langle t-\tau\rangle^{-\frac{1}{2}}\|\mathbf{u}(\tau)\cdot\nabla\theta(\tau)\|_{W^{10,1}}\dtau
\]
\[
+ \int_{0}^{t}\langle t-\tau\rangle^{-\frac{1}{2}}\|\partial_{1}(\mathbf{u}\cdot\nabla\omega_{2})(\tau)\|_{W^{6,1}}\dtau
+
\int_{0}^{t}\langle t-\tau\rangle^{-\frac{1}{2}}\|\partial_{2}(\mathbf{u}\cdot\nabla\omega_{1})(\tau)\|_{W^{6,1}}\dtau
\]
\[
+  \int_{0}^{t}\langle t-\tau\rangle^{-\frac{1}{2}}\|\partial_{1}(\boldsymbol{\omega}\cdot\nabla u_{2})(\tau)\|_{W^{6,1}}\dtau
+ \int_{0}^{t}\langle t-\tau\rangle^{-\frac{1}{2}}\|\partial_{2}(\boldsymbol{\omega}\cdot\nabla u_{1})(\tau)\|_{W^{6,1}}\dtau.
\]
Similar to the $L^\infty$-estimate of $\nabla_h \theta$ in Lemma \ref{ELD2},
\beq\label{EN2s}
\|\nabla_{h}\theta(t)\|_{H^{3}}
\ls \langle t\rangle^{-1}\left(\mathcal{E}_{0} + \mathcal{E}_{0}^{2}+\mathcal{E}^{2}(t)+\mathcal{E}^{\frac{12}{5}}(t)\right),
\eeq
\beq\label{EN2}
\|\theta(t)\|_{H^{5}}\ls \langle t\rangle^{-\frac{1}{2}}\left(\mathcal{E}_{0} + \mathcal{E}_{0}^{2}+\mathcal{E}^{2}(t)+\mathcal{E}^{\frac{12}{5}}(t)\right).
\eeq

To estimate $\omega_{1}$, we use (\ref{eLL3}), (\ref{EN2s}) and the fact $k\geq1$ to get
\[
\|\omega_{1}(t)\|_{H^{3}}
\lesssim \|e^{-\left(\xi^{2}+\eta^2+ \pi^2k ^{2}\right) t}\widehat{\omega_{10}}\|_{\widehat{H}^{3}}  + \int_{0}^{t}\|e^{-\left(\xi^{2}+\eta^2+ \pi^2 k^{2}\right)(t-\tau)}\widehat{\mathbf{u}\cdot\nabla\omega_{1}}(\tau)\|_{\widehat{H}^{3}} \dtau
\]
\[
 +
\int_{0}^{t}\|e^{-\left(\xi^{2}+\eta^2+\pi^2 k ^{2}\right)(t-\tau)}\widehat{\boldsymbol{\omega}\cdot\nabla u_{1}}(\tau)\|_{\widehat{H}^{3}} \dtau
+\int_{0}^{t}\|e^{-\left(\xi^{2}+\eta^2+ \pi^2k ^{2}\right)(t-\tau)}\widehat{\p_{2}\theta}(\tau)\|_{\widehat{H}^{3}}\dtau
\]
\[
\lesssim e^{- t}\|\boldsymbol{\omega}_{0}\|_{H^{3}} +
\int_{0}^{t}e^{-(t-\tau)}\|\mathbf{u}\cdot\nabla\omega_{1}(\tau)\|_{H^{3}}\dtau
\]
\[
+ \int_{0}^{t}e^{-(t-\tau)}\|\boldsymbol{\omega}\cdot\nabla u_{1}(\tau)\|_{H^{3}}\dtau+
\int_{0}^{t}e^{-(t-\tau)}\langle\tau\rangle^{-1}\langle\tau\rangle\|\p_{2}\theta(\tau)\|_{H^{3}}\dtau
\]
\[
\lesssim e^{- t}\|\boldsymbol{\omega}_{0}\|_{H^{3}} +
\int_{0}^{t}e^{-(t-\tau)}\|\mathbf{u}(\tau)\cdot\nabla \omega_{1}(\tau)\|_{W^{5,1}}\dtau + \int_{0}^{t}e^{-(t-\tau)}\|\boldsymbol{\omega}(\tau)\cdot\nabla u_{1}(\tau)\|_{W^{5,1}}\dtau
\]
\[
+ \langle t\rangle^{-1}
\left(\mathcal{E}_{0} + \mathcal{E}_{0}^{2}+\mathcal{E}^{2}(t)+\mathcal{E}^{\frac{12}{5}}(t)\right).
\]
Indeed, we obtain the bound
\[
\|\mathbf{u}\cdot\nabla \omega_{1}\|_{W^{5,1}}\ls \|\mathbf{u}\|_{H^{5}}\|\nabla\omega_{1}\|_{L^{2}}+\|\mathbf{u}\|_{L^{2}}\|\nabla\omega_{1}\|_{H^{5}}
\]
\[
\ls \|\mathbf{u}\|_{H^{5}}\|\nabla\omega_{1}\|_{L^{2}}+\|\mathbf{u}\|_{L^{2}}\|\nabla\omega_{1}\|_{L^{2}}^{\frac{1}{2}}
\|\nabla\omega_{1}\|_{H^{10}}^{\frac{1}{2}}.
\]
Then
\[
\int_{0}^{t}e^{-(t-\tau)}\|\mathbf{u}(\tau)\cdot\nabla \omega_{1}(\tau)\|_{W^{5,1}}\dtau\ls
\int_{0}^{t}e^{-(t-\tau)}\langle\tau\rangle^{-1}\dtau\mathcal{E}_{1}(t)\mathcal{E}_{3}(t)
\]
\[
+\int_{0}^{t}e^{-(t-\tau)}\langle\tau\rangle^{-1}\dtau\mathcal{E}_{1}^{\frac{1}{2}}(t)\mathcal{E}_{3}^{\frac{3}{2}}(t)
\ls \langle t\rangle^{-1}\mathcal{E}^{2}(t).
\]
Similarly,
\[
\int_{0}^{t}e^{-(t-\tau)}\|\boldsymbol{\omega}(\tau)\cdot\nabla u_{1}(\tau)\|_{W^{5,1}}\dtau\ls\langle t\rangle^{-1}\mathcal{E}^{2}(t).
\]
Hence
\beq\label{EN7}
\|\omega_{1}(t)\|_{H^{3}}\ls\langle t\rangle^{-1}
\left(\mathcal{E}_{0} + \mathcal{E}_{0}^{2}+\mathcal{E}^{2}(t)+\mathcal{E}^{\frac{12}{5}}(t)\right).
\eeq
Similar to the $H^3$-estimate of $\omega_1$,
\beq\label{EN8}
\|\omega_{2}(t)\|_{H^{3}}\ls\langle t\rangle^{-1}
\left(\mathcal{E}_{0} + \mathcal{E}_{0}^{2}+\mathcal{E}^{2}(t)+\mathcal{E}^{\frac{12}{5}}(t)\right).
\eeq
In a similar way,
\beq\label{EN7+}
\|\nabla_{h}\omega_{1}(t)\|_{H^{1}} + \|\nabla_{h}\omega_{2}(t)\|_{H^{1}} \ls\langle t\rangle^{-\frac{5}{4}}
\left(\mathcal{E}_{0} + \mathcal{E}_{0}^{2}+\mathcal{E}^{2}(t)+\mathcal{E}^{\frac{12}{5}}(t)\right).
\eeq
Furthermore, for
$(\xi,\eta)\in\mathbb{R}\times\mathbb{R}$ and $k\geq0$,
\[
\|\Lambda^{-1}\omega_{3}(t)\|_{L^{2}}
\ls \|e^{-\left(\xi^{2}+\eta^2+\pi^2k ^{2}\right) t}\left(\xi^{2}+\eta^2+ \pi^{2}k ^{2}\right)^{-\frac{1}{2}}\widehat{\omega_{30}}\|_{\widehat{L}^{2}}
\]
\[
+
\int_{0}^{t}\|e^{-\left(\xi^{2}+\eta^2+ \pi^{2}k ^{2}\right)(t-\tau)}\left(\xi^{2}+\eta^2+ \pi^{2}k ^{2}\right)^{-\frac{1}{2}}\widehat{\boldsymbol{\omega}\cdot\nabla u_{3}}(\tau)\|_{\widehat{L}^{2}}\dtau
\]
\[
 +\int_{0}^{t}\|e^{-\left(\xi^{2}+\eta^2+ \pi^{2}k ^{2}\right)(t-\tau)}\left(\xi^{2}+\eta^2+ \pi^{2}k ^{2}\right)^{-\frac{1}{2}}\widehat{\mathbf{u}\cdot\nabla\omega_{3}}(\tau)\|_{\widehat{L}^{2}}\dtau
\]
\beq\label{EN9}
:= P_{1}+P_{2}+P_{3}.
\eeq
Due to $k\geq 0$, we have
$P_1\ls \|\Lambda^{-1}\omega_{30}\|_{L^{2}}$ and
\[
P_1\ls \left(\int_{\mathbb{R}^{2}}e^{-2\left(\xi^{2}+\eta^{2}\right)t}\dxi\deta + \int_{\mathbb{R}^{3}}e^{-2\left(\xi^{2}+\eta^{2}+\zeta^{2}\right)t}\dxi\deta\dzeta \right)^{\frac{1}{2}}\|\widehat{\Lambda^{-1}\omega_{30}}\|_{\widehat{L}^{\infty}}
\]
\[
\ls  t^{-\frac{1}{2}}\|\Lambda^{-1}\omega_{30}\|_{L^{1}}.
\]
This implies that
\beq\label{EN13}
P_1 \ls  \langle t \rangle^{-\frac{1}{2}}\mathcal{E}_{0}.
\eeq
Note that
\[
P_2\ls
\int_{0}^{t}(t-\tau)^{-\frac{1}{2}}\|\widehat{\boldsymbol{\omega}u_{3}}(\tau)\|_{\widehat{L}^{\infty}}\dtau
\]
\[
\ls
\int_{0}^{t}(t-\tau)^{-\frac{1}{2}}\|\boldsymbol{\omega}(\tau)u_{3}(\tau)\|_{L^{1}}\dtau
\ls \int_{0}^{t}(t-\tau)^{-\frac{1}{2}}\|\boldsymbol{\omega}(\tau)\|_{L^{2}}\|u_{3}(\tau)\|_{L^{2}}\dtau.
\]
Also, we have
\[
P_2\ls
\int_{0}^{t}\|\widehat{\boldsymbol{\omega}u_{3}}(\tau)\|_{\widehat{L}^{2}}\dtau\ls
\int_{0}^{t}\|\boldsymbol{\omega}(\tau)u_{3}(\tau)\|_{L^{2}}\dtau\ls\int_{0}^{t}\|\boldsymbol{\omega}(\tau)\|_{H^{2}}\|u_{3}(\tau)\|_{H^{2}}\dtau.
\]
In consequence,
\beq\label{EN17}
P_2\ls
\int_{0}^{t}\langle t-\tau\rangle^{-\frac{1}{2}}\langle\tau\rangle^{-\frac{9}{4}}\dtau\mathcal{E}_{3}^{2}(t)
\ls \langle t\rangle^{-\frac{1}{2}} \mathcal{E}^{2}(t).
\eeq
Similarly,
\beq\label{EN20}
P_3\ls \langle t\rangle^{-\frac{1}{2}} \mathcal{E}^{2}(t).
\eeq
By substituting (\ref{EN13})--(\ref{EN20}) into (\ref{EN9}),
\beq\label{EN21}
\|\Lambda^{-1}\omega_{3}(t)\|_{L^{2}}\ls
\langle t\rangle^{-\frac{1}{2}} \left(\mathcal{E}_{0}+\mathcal{E}^{2}(t)\right).
\eeq
Note that
\[
\|\omega_{3}(t)\|_{L^{2}}
\ls \|e^{-\left(\xi^{2}+\eta^2+ \pi^{2}k ^{2}\right) t}\left(\xi^{2}+\eta^2+ \pi^{2}k ^{2}\right)^{\frac{1}{2}}\left(\xi^{2}+\eta^2+ \pi^{2}k ^{2}\right)^{-\frac{1}{2}}\widehat{\omega_{30}}\|_{\widehat{L}^{2}}
\]
\[
+
\int_{0}^{t}\|e^{-\left(\xi^{2}+\eta^2+ \pi^{2}k ^{2}\right)(t-\tau)}\left(\xi^{2}+\eta^2+ \pi^{2}k ^{2}\right)^{\frac{1}{2}}\left(\xi^{2}+\eta^2+ \pi^{2}k ^{2}\right)^{-\frac{1}{2}}\widehat{\boldsymbol{\omega}\cdot\nabla u_{3}}(\tau)\|_{\widehat{L}^{2}}\dtau
\]
\[
 +\int_{0}^{t}\|e^{-\left(\xi^{2}+\eta^2+ \pi^{2}k ^{2}\right)(t-\tau)}\left(\xi^{2}+\eta^2+ \pi^{2}k ^{2}\right)^{\frac{1}{2}}\left(\xi^{2}+\eta^2+ \pi^{2}k ^{2}\right)^{-\frac{1}{2}}\widehat{\mathbf{u}\cdot\nabla\omega_{3}}(\tau)\|_{\widehat{L}^{2}}\dtau,
\]
where
$(\xi,\eta)\in\mathbb{R}\times\mathbb{R}$ and $k\geq0$. Thus,
\[
\|\omega_{3}(t)\|_{L^{2}}
\ls
\langle t\rangle^{-1}\|\Lambda^{-1}\omega_{30}\|_{L^{1}}
\]
\[
+
\int_{0}^{t}\langle t-\tau\rangle^{-1}\|\boldsymbol{\omega}(\tau)\|_{H^{2}}\|u_{3}(\tau)\|_{H^{2}}\dtau
+\int_{0}^{t}\langle t-\tau\rangle^{-1}\|\mathbf{u}(\tau)\|_{H^{2}}\|\omega_{3}(\tau)\|_{H^{2}}\dtau
\]
\[
\ls
\langle t\rangle^{-1}\mathcal{E}_{0}+
\int_{0}^{t}\langle t-\tau\rangle^{-1}\langle\tau\rangle^{-\frac{9}{4}}\dtau\mathcal{E}_{3}^{2}(t)
+\int_{0}^{t}\langle t-\tau\rangle^{-1}\langle\tau\rangle^{-\frac{3}{2}}\dtau\mathcal{E}_{3}^{2}(t)
\]
\beq\label{EX1}
\ls \langle t\rangle^{-1}\left(\mathcal{E}_{0}+\mathcal{E}^{2}(t)\right).
\eeq
In the similar way,
\[
\|\omega_{3}(t)\|_{\dot{H}^{1}}\ls
\langle t\rangle^{-\frac{3}{2}}\mathcal{E}_{0}
\]
\[
+
\int_{0}^{t}\langle t-\tau\rangle^{-\frac{3}{2}}\langle\tau\rangle^{-\frac{9}{4}}\dtau\mathcal{E}_{3}^{2}(t)
+\int_{0}^{t}\langle t-\tau\rangle^{-\frac{3}{2}}\langle\tau\rangle^{-\frac{3}{2}}\dtau\mathcal{E}_{3}^{2}(t)
\]
\beq\label{EX2}
\ls \langle t\rangle^{-\frac{3}{2}}\left(\mathcal{E}_{0}+\mathcal{E}^{2}(t)\right),
\eeq
\[
\|\omega_{3}(t)\|_{\dot{H}^{2}}\ls
\langle t\rangle^{-2}\mathcal{E}_{0}
\]
\[
+
\int_{0}^{t}\langle t-\tau\rangle^{-2}\langle\tau\rangle^{-\frac{9}{4}}\dtau\mathcal{E}_{3}^{2}(t)
+\int_{0}^{t}\langle t-\tau\rangle^{-2}\langle\tau\rangle^{-\frac{3}{2}}\dtau\mathcal{E}_{3}^{2}(t)
\]
\beq\label{EX3}
\ls \langle t\rangle^{-\frac{3}{2}}\left(\mathcal{E}_{0}+\mathcal{E}^{2}(t)\right).
\eeq
Furthermore,
\[
\|\omega_{3}(t)\|_{\dot{H}^{3}}\ls
\langle t\rangle^{-\frac{5}{2}}\|\Lambda^{-1}\omega_{30}\|_{L^{1}}
+
\int_{0}^{t}\langle t-\tau\rangle^{-\frac{5}{2}}\|\boldsymbol{\omega}(\tau)\|_{H^{3}}\|u_{3}(\tau)\|_{H^{4}}\dtau
\]
\[
+\int_{0}^{t}\langle t-\tau\rangle^{-\frac{5}{2}}\|\mathbf{u}(\tau)\|_{H^{3}}\|\omega_{3}(\tau)\|_{H^{4}}\dtau
\]
\[
\ls \langle t\rangle^{-\frac{5}{2}}\mathcal{E}_{0}+
\int_{0}^{t}\langle t-\tau\rangle^{-\frac{5}{2}}\langle\tau\rangle^{-\frac{13}{8}}\dtau\mathcal{E}_{1}^{\frac{1}{2}}(t)\mathcal{E}_{3}^{\frac{3}{2}}(t)
\]
\beq\label{EX4}
+\int_{0}^{t}\langle t-\tau\rangle^{-\frac{5}{2}}\langle\tau\rangle^{-1}\dtau
\mathcal{E}_{1}^{\frac{1}{2}}(t)\mathcal{E}_{3}^{\frac{3}{2}}(t)
\ls \langle t\rangle^{-1}\left(\mathcal{E}_{0}+\mathcal{E}^{2}(t)\right),
\eeq
where we use
\[
\|u_{3}(\tau)\|_{H^{4}}\ls\|u_{3}(\tau)\|_{H^{3}}^{\frac{1}{2}}\|u_{3}(\tau)\|_{H^{5}}^{\frac{1}{2}},\,\|\omega_{3}(\tau)\|_{H^{4}}\ls\|\omega_{3}(\tau)\|_{H^{3}}^{\frac{1}{2}}\|\omega_{3}(\tau)\|_{H^{5}}^{\frac{1}{2}}.
\]
Therefore, one infers from (\ref{EX1})--(\ref{EX4}) that
\beq\label{EX5}
\|\omega_3(t)\|_{H^{3}}\ls \langle t\rangle^{-1}\left(\mathcal{E}_{0}+\mathcal{E}^{2}(t)\right).
\eeq
Due to (\ref{Bou3})$_{3}$, we have
\[
\|u_{1}(t)\|_{H^{4}}\ls\|\mathcal{F}\left(\partial_{2}(-\Delta)^{-1}\omega_{3}\right)(t)\|_{\widehat{H}^{4}}+\|\mathcal{F}\left(\partial_{3}(-\Delta)^{-1}\omega_{2}\right)(t)\|_{\widehat{H}^{4}}.
\]
Moreover,
\[
\|\mathcal{F}\left(\partial_{2}(-\Delta)^{-1}\omega_{3}\right)(t)\|_{\widehat{H}^{4}}=
\left\|\frac{i\eta}{\xi^2+\eta^2+\pi^{2} k^2}\widehat{\omega_{3}}(t)\right\|_{\widehat{H}^{4}}
\]
\[
\ls \|\widehat{\Lambda^{-1}\omega_{3}}(t)\|_{\widehat{L}^{2}} + \|\widehat{\omega_{3}}(t)\|_{\widehat{H}^{3}},\,(\xi,\eta)\in\mathbb{R}
\times\mathbb{R},\,k\geq0,
\]
\[
\|\mathcal{F}\left(\partial_{3}(-\Delta)^{-1}\omega_{2}\right)(t)\|_{\widehat{H}^{4}}
=
\left\|\frac{k}{\xi^2+\eta^2+\pi^{2} k^2}\widehat{\omega_{2}}(t)\right\|_{\widehat{H}^{4}}
\]
\[
\ls \|\widehat{\omega_{2}}(t)\|_{\widehat{H}^{3}},\,(\xi,\eta)\in\mathbb{R}
\times\mathbb{R},\,k\geq 1.
\]
Then combining (\ref{EN8}), (\ref{EN21}) and (\ref{EX5}),
\[
\|u_{1}(t)\|_{H^{4}}
\ls
\langle t\rangle^{-\frac{1}{2}}
\left(\mathcal{E}_{0} + \mathcal{E}_{0}^{2}+\mathcal{E}^{2}(t)+\mathcal{E}^{\frac{12}{5}}(t)\right).
\]
In the similar way,
\[
\|u_{2}(t)\|_{H^{4}}
\ls
\langle t\rangle^{-\frac{1}{2}}
\left(\mathcal{E}_{0} + \mathcal{E}_{0}^{2}+\mathcal{E}^{2}(t)+\mathcal{E}^{\frac{12}{5}}(t)\right).
\]
Note that
\[
\|u_{3}(t)\|_{H^{3}}\ls
\|\mathcal{F}\left(\partial_{1}(-\Delta)^{-1}\omega_{2}\right)(t)\|_{\widehat{H}^{3}}
+ \|\mathcal{F}\left(\partial_{2}(-\Delta)^{-1}\omega_{1}\right)(t)\|_{\widehat{H}^{3}}
\]
\[
\ls\|\widehat{\partial_{1}\omega_{2}}(t)\|_{\widehat{H}^{1}}+\|\widehat{\partial_{2}\omega_{1}}(t)\|_{\widehat{H}^{1}}.
\]
By using (\ref{EN7+}),
\[
\|u_{3}(t)\|_{H^{3}}\ls
\langle t\rangle^{-\frac{5}{4}}
\left(\mathcal{E}_{0} + \mathcal{E}_{0}^{2}+\mathcal{E}^{2}(t)+\mathcal{E}^{\frac{12}{5}}(t)\right).
\]
Thus we conclude the proof of Lemma \ref{ELD3}.
\hfill$\square$
\begin{Lemma}\label{ELD4}
Let $m\geq 31$. Then
\begin{equation}\label{ET1}
\mathcal{E}_{4}(t)\ls \mathcal{E}_{0} + \mathcal{E}_{0}^{2}+\mathcal{E}^{2}(t)+\mathcal{E}^{\frac{12}{5}}(t).
\end{equation}
\end{Lemma}
\proof
Using the fact that $\mathcal{L}_{2}(0) = 0$, one deduces from (\ref{if5}) that
\[
\partial_{t}\theta(t,\mathbf{x})= \partial_t\mathcal{L}_{1}(t)\theta_{0}(\mathbf{x}) + \partial_t\mathcal{L}_{2}(t)\left(\frac{1}{2}(-\Delta)\theta_{0}(\mathbf{x}) +\p_{t}\theta(0,\mathbf{x})\right)
+ \int_{0}^{t}\partial_{\tau}\mathcal{L}_{2}(t-\tau)f_{2}(\tau,\mathbf{x})\dtau.
\]
By applying Lemma \ref{AL2},
\[
\|\partial_{t}\theta(t)\|_{H^{4}}
\ls\langle t\rangle^{-\frac{3}{2}}\|\boldsymbol{\omega}_{0}\|_{W^{7,1}} + \langle t\rangle^{-\frac{3}{2}}\|\theta_{0}\|_{W^{10,1}}
\]
\[
 +
\langle t\rangle^{-\frac{3}{2}}\|\theta_{0}\|_{H^{m}}
(\|\boldsymbol{\omega}_{0}\|_{H^{m}}+\|\Lambda^{-1}\omega_{30}\|_{L^{2}})
+ \int_{0}^{t}\langle t-\tau \rangle^{-\frac{3}{2}}\|\partial_{\tau}\mathbf{u}(\tau)\cdot\nabla\theta(\tau)\|_{W^{8,1}}\dtau
\]
\[
+  \int_{0}^{t}\langle t-\tau\rangle^{-\frac{3}{2}}\|\mathbf{u}(\tau)\cdot\nabla\partial_{\tau}\theta(\tau)\|_{W^{8,1}}\dtau
+  \int_{0}^{t}\langle t-\tau\rangle^{-\frac{3}{2}}\|\mathbf{u}(\tau)\cdot\nabla\theta(\tau)\|_{W^{10,1}}\dtau
\]
\[
+ \int_{0}^{t}\langle t-\tau\rangle^{-\frac{3}{2}}\|\partial_{1}(\mathbf{u}\cdot\nabla\omega_{2})(\tau)\|_{W^{6,1}}\dtau
+
\int_{0}^{t}\langle t-\tau\rangle^{-\frac{3}{2}}\|\partial_{2}(\mathbf{u}\cdot\nabla\omega_{1})(\tau)\|_{W^{6,1}}\dtau
\]
\[
+  \int_{0}^{t}\langle t-\tau\rangle^{-\frac{3}{2}}\|\partial_{1}(\boldsymbol{\omega}\cdot\nabla u_{2})(\tau)\|_{W^{6,1}}\dtau
+ \int_{0}^{t}\langle t-\tau\rangle^{-\frac{3}{2}}\|\partial_{2}(\boldsymbol{\omega}\cdot\nabla u_{1})(\tau)\|_{W^{6,1}}\dtau.
\]
Similar to the $L^\infty$-estimate of $\nabla_h \theta$ in Lemma \ref{ELD2},
\beq\label{ET4}
\|\partial_{t}\theta(t)\|_{H^{4}}
\ls \langle t\rangle^{-\frac{5}{4}}\left(\mathcal{E}_{0} + \mathcal{E}_{0}^{2}+\mathcal{E}^{2}(t)+\mathcal{E}^{\frac{12}{5}}(t)\right).
\eeq

Taking time derivative on (\ref{if1})$_1$ yields
\[
\partial_{tt}\omega_1-\Delta\partial_{t}\omega_1 = \partial_{2}\partial_{t}\theta + \partial_{t}\left(\boldsymbol{\omega}\cdot \nabla u_1\right)- \partial_{t}\left(\mathbf{u}\cdot \nabla \omega_1\right),
\]
\[
\partial_{tt}\omega_2-\Delta\partial_{t}\omega_2 = -\partial_{1}\partial_{t}\theta + \partial_{t}\left(\boldsymbol{\omega}\cdot \nabla u_2\right)- \partial_{t}\left(\mathbf{u}\cdot \nabla \omega_2\right).
\]
By the Duhamel's principle,
\[
   \partial_{t}\omega_{1}(t,\mathbf{x})= e^{t\Delta }\partial_{t}\omega_{1}(0,\mathbf{x}) + \int_{0}^{t}e^{(t-\tau)\Delta}\left(\partial_{2}\partial_{\tau}\theta +\partial_{\tau}(\boldsymbol{\omega}\cdot \nabla u_{1})- \partial_{\tau}(\mathbf{u}\cdot \nabla \omega_{1})\right)(\tau,\mathbf{x})\dtau,
\]
\[
   \partial_{t}\omega_{2}(t,\mathbf{x})= e^{t\Delta }\partial_{t}\omega_{2}(0,\mathbf{x}) + \int_{0}^{t}e^{(t-\tau)\Delta}\left(-\partial_{1}\partial_{\tau}\theta+ \partial_{\tau}(\boldsymbol{\omega}\cdot \nabla u_{2}) - \partial_{\tau}(\mathbf{u}\cdot \nabla \omega_{2})\right)(\tau,\mathbf{x})\dtau.
\]
From (\ref{if1})$_1$ we find that
\beq\label{EM2}
\partial_{t}\boldsymbol{\omega}(0,\mathbf{x})=\Delta\boldsymbol{\omega}_{0}(\mathbf{x})-\mathbf{u}_{0}\cdot \nabla \boldsymbol{\omega}_{0}(\mathbf{x})+\boldsymbol{\omega}_{0}\cdot \nabla \mathbf{u}_{0}(\mathbf{x}) + \left(\partial_{2}\theta_{0}(\mathbf{x}), -\partial_{1}\theta_{0}(\mathbf{x}), 0\right).
\eeq
Moreover, we obtain
\[
\|\partial_{t}\omega_{1}(t)\|_{L^{2}}
\lesssim \|e^{-\left(\xi^{2}+\eta^2+ \pi^2k ^{2}\right) t}\widehat{\partial_{t}\omega_{1}}(0)\|_{\widehat{L}^{2}}
+\int_{0}^{t}\|e^{-\left(\xi^{2}+\eta^2+\pi^2 k ^{2}\right)(t-\tau)}\widehat{\p_\tau(\boldsymbol{\omega}\cdot\nabla u_{1})}(\tau)\|_{\widehat{L}^{2}} \dtau
\]
\[
+\int_{0}^{t}\|e^{-\left(\xi^{2}+\eta^2+ \pi^2k ^{2}\right)(t-\tau)}\widehat{\p_\tau(\mathbf{u}\cdot\nabla\omega_{1})}(\tau)\|_{\widehat{L}^{2}} \dtau
+\int_{0}^{t}\|e^{-\left(\xi^{2}+\eta^2+ \pi^2k ^{2}\right)(t-\tau)}\widehat{\p_{2}\p_\tau\theta}(\tau)\|_{\widehat{L}^{2}}\dtau
\]
\[
\lesssim e^{- t}\|\p_t\omega_{1}(0)\|_{L^{2}}+
\int_{0}^{t}e^{-(t-\tau)}\|\partial_{\tau}(\boldsymbol{\omega}\cdot \nabla u_{1})(\tau)\|_{L^{2}}\dtau
\]
\[
+ \int_{0}^{t}e^{-(t-\tau)}\|\partial_{\tau}(\mathbf{u}\cdot \nabla \omega_{1})(\tau)\|_{L^{2}}\dtau+
\int_{0}^{t}e^{-(t-\tau)}\|\p_{2}\p_\tau\theta(\tau)\|_{L^{2}}\dtau
\]
\beq\label{EM1}
:= M_1+M_2+M_3+M_4.
\eeq
By using (\ref{EM2}),
\beq\label{EM9}
M_1\ls e^{- t}(\mathcal{E}_{0}+\mathcal{E}_{0}^{2}).
\eeq
For $M_{2}$, we calculate that
\[
\|\partial_{t}(\boldsymbol{\omega}\cdot \nabla u_{1})\|_{L^{2}}\ls \|\partial_{t}\boldsymbol{\omega}\cdot \nabla u_{1}\|_{L^{2}}+\|\boldsymbol{\omega}\cdot \nabla\partial_{t}u_{1}\|_{L^{2}}
\]
\[
\ls \|\partial_{t}\boldsymbol{\omega}\|_{L^{2}}\|\nabla u_{1}\|_{H^{2}}+\|\boldsymbol{\omega}\|_{H^{2}}\|\nabla\partial_{t}u_{1}\|_{L^{2}}.
\]
Then
\[
M_2\ls \int_{0}^{t}e^{-(t-\tau)}\langle\tau\rangle^{-\frac{7}{4}}
\dtau\mathcal{E}_{4}(t)\mathcal{E}_{3}(t)
+\int_{0}^{t}e^{-(t-\tau)}
\langle\tau\rangle^{-\frac{9}{4}}\dtau\mathcal{E}_{4}(t)\mathcal{E}_{3}(t)
\]
\beq\label{EM3}
\ls \left(\langle t\rangle^{-\frac{7}{4}}+\langle t\rangle^{-\frac{9}{4}}\right)\mathcal{E}^2(t)\ls\langle t\rangle^{-\frac{5}{4}}\mathcal{E}^2(t).
\eeq
Similarly,
\[
M_3
\ls \int_{0}^{t}e^{-(t-\tau)}\|\partial_{\tau}\mathbf{u}(\tau)\|_{L^{2}}
\|\omega_{1}(\tau)\|_{H^{3}}\dtau
+\int_{0}^{t}e^{-(t-\tau)}\|\mathbf{u}(\tau)\|_{L^{2}}
\|\p_{\tau}\omega_{1}(\tau)\|_{H^{3}}\dtau
\]
\[
\ls \int_{0}^{t}e^{-(t-\tau)}\langle\tau\rangle^{-\frac{9}{4}}\dtau\mathcal{E}_{4}(t)\mathcal{E}_{3}(t)
+\int_{0}^{t}e^{-(t-\tau)}
\langle\tau\rangle^{-\frac{5}{4}}\dtau\mathcal{E}_{3}(t)\mathcal{E}_{4}^{\frac{3}{5}}(t)\mathcal{E}_{1}^{\frac{2}{5}}(t)
\]
\beq\label{EM4}
 \ls \left(\langle t\rangle^{-\frac{9}{4}}+\langle t\rangle^{-\frac{5}{4}}\right)\mathcal{E}^2(t)\ls\langle t\rangle^{-\frac{5}{4}}\mathcal{E}^2(t),
\eeq
where we use
\[
\| \p_{\tau}\omega_{1}(\tau)\|_{H^{3}}\ls\| \p_{\tau}\omega_{1}(\tau)\|_{L^{2}}^{\frac{3}{5}}\| \p_{\tau}\omega_{1}(\tau)\|_{H^{\frac{15}{2}}}^{\frac{2}{5}}.
\]
For $M_{4}$, estimate (\ref{ET4}) implies that
\[
M_4\ls \int_{0}^{t}e^{-(t-\tau)}\langle \tau \rangle^{-\frac{5}{4}}\langle \tau \rangle^{\frac{5}{4}}\|\p_{2}\p_\tau\theta(\tau)\|_{L^{2}}\dtau
\]
\beq\label{EM5}
\ls \langle t\rangle^{-\frac{5}{4}}\left(\mathcal{E}_{0} + \mathcal{E}_{0}^{2}+\mathcal{E}^{2}(t)+\mathcal{E}^{\frac{12}{5}}(t)\right).
\eeq
Inserting (\ref{EM9})--(\ref{EM5}) into (\ref{EM1}) yields
\beq\label{EM8}
\|\partial_{t}\omega_{1}(t)\|_{L^{2}}
\ls \langle t\rangle^{-\frac{5}{4}}\left(\mathcal{E}_{0} + \mathcal{E}_{0}^{2}+\mathcal{E}^{2}(t)+\mathcal{E}^{\frac{12}{5}}(t)\right).
\eeq
Similarly,
\beq\label{EM10}
\|\partial_{t}\omega_{2}(t)\|_{L^{2}}
\ls \langle t\rangle^{-\frac{5}{4}}\left(\mathcal{E}_{0} + \mathcal{E}_{0}^{2}+\mathcal{E}^{2}(t)+\mathcal{E}^{\frac{12}{5}}(t)\right).
\eeq
From (\ref{if4}), we see that
\[
   \omega_{3}(t,\mathbf{x})= e^{t\Delta }\omega_{30}(\mathbf{x}) + \int_{0}^{t}e^{(t-\tau)\Delta}\left((\boldsymbol{\omega}\cdot \nabla u_{3})(\tau, \mathbf{x})-(\mathbf{u}\cdot \nabla \omega_{3})(\tau, \mathbf{x})\right)\dtau.
\]
Then
\[
   \p_t\omega_{3}(t,\mathbf{x})=\Delta e^{t\Delta }\omega_{30}(\mathbf{x}) + \int_{0}^{t}\Delta e^{(t-\tau)\Delta}\left((\boldsymbol{\omega}\cdot \nabla u_{3})(\tau, \mathbf{x})-(\mathbf{u}\cdot \nabla \omega_{3})(\tau, \mathbf{x})\right)\dtau
\]
\[
   + (\boldsymbol{\omega}\cdot \nabla u_{3})(t, \mathbf{x})-(\mathbf{u}\cdot \nabla \omega_{3})(t, \mathbf{x}).
\]
For $(\xi,\eta)\in\mathbb{R}\times\mathbb{R}$ and $k\geq 0$,
\[
\|\partial_{t}\Lambda^{-1}\omega_{3}(t)\|_{L^{2}}
\lesssim \|e^{-\left(\xi^{2}+\eta^2+ \pi^2k ^{2}\right) t}\left(\xi^{2}+\eta^2+ \pi^2k ^{2}\right)\widehat{\Lambda^{-1}\omega_{30}}\|_{\widehat{L}^{2}}
\]
\[
+\int_{0}^{t}\|e^{-\left(\xi^{2}+\eta^2+\pi^2 k ^{2}\right)(t-\tau)}\left(\xi^{2}+\eta^2+ \pi^2k ^{2}\right)\left(\xi^{2}+\eta^2+\pi^2 k ^{2}\right)^{-\frac{1}{2}}\widehat{\boldsymbol{\omega}\cdot\nabla u_{3}}(\tau)\|_{\widehat{L}^{2}} \dtau
\]
\[
+\int_{0}^{t}\|e^{-\left(\xi^{2}+\eta^2+ \pi^2k ^{2}\right)(t-\tau)}\left(\xi^{2}+\eta^2+ \pi^2k ^{2}\right)\left(\xi^{2}+\eta^2+\pi^2 k ^{2}\right)^{-\frac{1}{2}}\widehat{\mathbf{u}\cdot\nabla\omega_{3}}(\tau)\|_{\widehat{L}^{2}} \dtau
\]
\[
+\|\left(\xi^{2}+\eta^2+\pi^2 k ^{2}\right)^{-\frac{1}{2}}\widehat{\boldsymbol{\omega}\cdot\nabla u_{3}}(t)\|_{\widehat{L}^{2}} + \|\left(\xi^{2}+\eta^2+\pi^2 k ^{2}\right)^{-\frac{1}{2}}\widehat{\mathbf{u}\cdot\nabla\omega_{3}}(t)\|_{\widehat{L}^{2}}
\]
\beq\label{EU1}
:=Y_1+ Y_2+ Y_3+ Y_4 + Y_5.
\eeq
For $Y_1$, one gets $Y_1\ls \|\omega_{30}\|_{H^{1}}$. Since $k\geq0$, we also obtain
\[
Y_1\ls \|e^{-\left(\xi^{2}+\eta^2+ \pi^2k ^{2}\right) t}\left(\xi^{2}+\eta^2+ \pi^2k ^{2}\right)\widehat{\Lambda^{-1}\omega_{30}}\|_{\widehat{L}^{2}}
\]
\[
\ls t^{-1}\left(\int_{\mathbb{R}^{2}}e^{-2\left(\xi^{2}+\eta^{2}\right)t}\dxi\deta + \int_{\mathbb{R}^{3}}e^{-2\left(\xi^{2}+\eta^{2}+\zeta^{2}\right)t}\dxi\deta\dzeta \right)^{\frac{1}{2}}\|\widehat{\Lambda^{-1}\omega_{30}}\|_{\widehat{L}^{\infty}}
\]
\[
\ls  t ^{-\frac{3}{2}}\|\Lambda^{-1}\omega_{30}\|_{L^{1}}.
\]
This implies that
\beq\label{YE1}
Y_1\ls  \langle t \rangle ^{-\frac{3}{2}}\mathcal{E}_{0}.
\eeq
On one hand, we have
\[
Y_2\ls \int_{0}^{t}\|\boldsymbol{\omega}(\tau)u_{3}(\tau)\|_{H^{2}} \dtau.
\]
On the other hand,
\[
Y_2\ls\int_{0}^{t}\langle t- \tau \rangle^{-\frac{3}{2}}\|\boldsymbol{\omega}(\tau)u_{3}(\tau)\|_{L^{1}} \dtau \ls
\int_{0}^{t}\langle t-\tau \rangle^{-\frac{3}{2}}\|\boldsymbol{\omega}(\tau)\|_{L^2}\| u_{3}(\tau)\|_{L^{2}} \dtau.
\]
Then
\[
Y_2\ls\int_{0}^{t}\langle t-\tau \rangle^{-\frac{3}{2}}\|\boldsymbol{\omega}(\tau)\|_{H^{2}}\| u_{3}(\tau)\|_{H^{2}} \dtau
\]
\beq\label{YE4}
\ls \int_{0}^{t}\langle t-\tau \rangle^{-\frac{3}{2}}\langle \tau \rangle^{-\frac{9}{4}}\dtau\mathcal{E}_{3}^{2}(t)
\ls
\langle t\rangle^{-\frac{3}{2}}\mathcal{E}^{2}(t).
\eeq
In the similar way,
\beq\label{YE5}
Y_3\ls \langle t\rangle^{-\frac{3}{2}}\mathcal{E}^{2}(t).
\eeq
For $Y_{4}$, we obtain the bound $Y_4 \ls \|\boldsymbol{\omega}u_{3}(t)\|_{L^{2}}\ls \mathcal{E}^{2}(t)$. Then
\beq\label{ec1}
Y_{4} \ls \langle t\rangle^{-1}\langle t\rangle\|\boldsymbol{\omega}(t)\|_{L^{2}}\langle t\rangle^{-\frac{5}{4}}\langle t\rangle^{\frac{5}{4}}\|u_{3}(t)\|_{H^{2}}
\ls \langle t\rangle^{-\frac{9}{4}}\mathcal{E}^{2}(t).
\eeq
Similarly for $Y_5$,
\beq\label{ec2}
Y_5 \ls \langle t\rangle^{-\frac{3}{2}}\mathcal{E}^{2}(t).
\eeq
Plugging (\ref{YE1})--(\ref{ec2}) into (\ref{EU1}) yields
\beq\label{ec3}
\|\partial_{t}\Lambda^{-1}\omega_{3}(t)\|_{L^{2}}
\ls \langle t\rangle^{-\frac{3}{2}}\left(\mathcal{E}_{0}^{2}+\mathcal{E}^{2}(t)\right).
\eeq
Now we turn to estimate $L^2$-estimate of $\p_t\omega_3$. For $(\xi,\eta)\in\mathbb{R}\times\mathbb{R}$ and $k\geq0$,
\[
\|\partial_{t}\omega_{3}(t)\|_{L^{2}}
\lesssim \|e^{-\left(\xi^{2}+\eta^2+ \pi^2k ^{2}\right) t}\left(\xi^{2}+\eta^2+ \pi^2k ^{2}\right)^{\frac{3}{2}}\left(\xi^{2}+\eta^2+ \pi^2k ^{2}\right)^{-\frac{1}{2}}\widehat{\omega_{30}}\|_{\widehat{L}^{2}}
\]
\[
+\int_{0}^{t}\|e^{-\left(\xi^{2}+\eta^2+\pi^2 k ^{2}\right)(t-\tau)}\left(\xi^{2}+\eta^2+ \pi^2k ^{2}\right)^{\frac{3}{2}}\left(\xi^{2}+\eta^2+ \pi^2k ^{2}\right)^{-\frac{1}{2}}\widehat{\boldsymbol{\omega}\cdot\nabla u_{3}}(\tau)\|_{\widehat{L}^{2}} \dtau
\]
\[
+\int_{0}^{t}\|e^{-\left(\xi^{2}+\eta^2+ \pi^2k ^{2}\right)(t-\tau)}\left(\xi^{2}+\eta^2+ \pi^2k ^{2}\right)^{\frac{3}{2}}\left(\xi^{2}+\eta^2+ \pi^2k ^{2}\right)^{-\frac{1}{2}}\widehat{\mathbf{u}\cdot\nabla\omega_{3}}(\tau)\|_{\widehat{L}^{2}} \dtau
\]
\[
+\|\widehat{\boldsymbol{\omega}\cdot\nabla u_{3}}(t)\|_{\widehat{L}^{2}} + \|\widehat{\mathbf{u}\cdot\nabla\omega_{3}}(t)\|_{\widehat{L}^{2}}.
\]
Then
\[
\|\partial_{t}\omega_{3}(t)\|_{L^{2}}
\lesssim \langle t\rangle^{-2}\mathcal{E}_{0}
+\int_{0}^{t}\langle t -\tau\rangle^{-2}\|\boldsymbol{\omega}(\tau)\|_{H^3}\|u_{3}(\tau)\|_{H^{3}} \dtau
\]
\[
+\int_{0}^{t}\langle t -\tau\rangle^{-2}\|\mathbf{u}(\tau)\|_{H^3}\|\omega_{3}(\tau)\|_{H^{3}}\dtau
+\|\boldsymbol{\omega}(t)\|_{H^{2}}\|\nabla u_{3}(t)\|_{L^{2}} + \|\mathbf{u}(t)\|_{H^{2}}\|\nabla\omega_{3}(t)\|_{L^{2}}
\]
\[
\ls \langle t\rangle^{-2}\mathcal{E}_{0}++\int_{0}^{t}\langle t -\tau\rangle^{-2}\langle \tau \rangle^{-\frac{9}{4}}\dtau\mathcal{E}_{3}^{2}(t)
+\int_{0}^{t}\langle t -\tau\rangle^{-2}\langle \tau \rangle^{-\frac{3}{2}}\dtau\mathcal{E}_{3}^{2}(t)
\]
\beq\label{we1}
+ \langle t\rangle^{-\frac{9}{4}}\mathcal{E}_{3}^{2}(t)+ \langle t\rangle^{-\frac{3}{2}}\mathcal{E}_{3}^{2}(t)
\ls \langle t\rangle^{-\frac{3}{2}}\left(\mathcal{E}_{0}+\mathcal{E}^{2}(t)\right).
\eeq
By using (\ref{EM10}) and (\ref{ec3})--(\ref{we1}),
\[
\|\p_t u_1\|_{H^{1}} \ls \|\p_t \omega_2\|_{L^{2}} + \|\p_t \Lambda^{-1}\omega_3\|_{L^{2}} + \|\p_t \omega_3\|_{L^{2}}
\ls \langle t\rangle^{-\frac{5}{4}}\left(\mathcal{E}_{0} + \mathcal{E}_{0}^{2}+\mathcal{E}^{2}(t)+\mathcal{E}^{\frac{12}{5}}(t)\right).
\]
In the similar way,
\[
\|\p_t u_2\|_{H^{1}}
\ls \langle t\rangle^{-\frac{5}{4}}\left(\mathcal{E}_{0} + \mathcal{E}_{0}^{2}+\mathcal{E}^{2}(t)+\mathcal{E}^{\frac{12}{5}}(t)\right).
\]
Furthermore, (\ref{EM8}) and (\ref{EM10}) imply that
\[
\|\p_t u_3\|_{H^{1}} \ls \|\p_t \omega_2\|_{L^{2}} + \|\p_t \omega_1\|_{L^{2}}
\ls \langle t\rangle^{-\frac{5}{4}}\left(\mathcal{E}_{0} + \mathcal{E}_{0}^{2}+\mathcal{E}^{2}(t)+\mathcal{E}^{\frac{12}{5}}(t)\right).
\]
Thus the proof of Lemma \ref{ELD4} is completed.
\hfill$\square$
\subsection{Proof of Theorem \ref{Mrt1}}
We now apply the bootstrap argument to conclude the proof of Theorem \ref{Mrt1}. According to Lemma \ref{ELD1}--\ref{ELD4}, there exists a constant $C_{0}\ge 1$ such that
\begin{equation}\label{Ep362+}
\mathcal{E}(t) \leq C_{0}\left(\mathcal{E}_{0}+\mathcal{E}_{0}^{2}\right)+ C_{0}\left(\mathcal{E}^{2}(t) +\mathcal{E}^{\frac{12}{5}}(t) \right).
\end{equation}
Assume that
\begin{equation}\label{assim2}
\|\theta_{0}\|_{W^{10,1}} + \|\theta_{0}\|_{H^{m+1}} +\|\boldsymbol{\omega}_{0}\|_{W^{7,1}} + \|\boldsymbol{\omega}_{0}\|_{H^{m}}+\|\Lambda^{-1}\omega_{3_{0}}\|_{L^{2}} + \|\Lambda^{-1}\omega_{3_{0}}\|_{L^{1}}\leq \epsilon_{0}
\end{equation}
with $\epsilon_{0}\in(0,1)$  to be determined later.
Then by the definition of $E_i(t), i=1,2,3,4$, there exists $C_1\ge 2C_0$ such that
\[
\mathcal{E}(0)= \mathcal{E}_{1}(0)+\mathcal{E}_{2}(0)+\mathcal{E}_{3}(0)+\mathcal{E}_{4}(0)
\leq C_1\mathcal{E}_{0}\leq C_{1}\epsilon_{0}.
\]
From Proposition \ref{Lth1}, there exists $T\in (0,T^{*})$ such that
\[
(\boldsymbol{\omega}, \Lambda^{-1}\omega_{3}, \theta) \in C([0,T]; H^{m})\times C([0,T];L^{2})\times C([0,T];H^{m+1}).
\]
Now we assume that
\begin{equation}\label{boots1}
\mathcal{E}(t) \leq 4C_{1}\epsilon_{0}\leq1,\,\text{ for all }t\in[0,T].
\end{equation}
By using (\ref{Ep362+})--(\ref{assim2}) and the bootstrap hypothesis (\ref{boots1}),
\[
\mathcal{E}(t)\leq 2C_{0}\mathcal{E}_{0}+ 2C_{0}\mathcal{E}^{2}(t)
\leq C_{1}\epsilon_{0}+  32 C_0 C_{1}^{2}\epsilon_{0}^{2},\,\text{ for all}~t \in [0,T].
\]
Then taking $\epsilon_{0}$ so small that $32 C_0 C_{1}\epsilon_{0}\leq 1$,
\[
\mathcal{E}(t)\leq 2C_{1} \epsilon_{0},\,\text{ for all}~t \in [0,T].
\]
Thus the bootstrap is closed. By repeating the argument, the local solution is prolonged to the whole time interval $[0,\infty)$ and all decay estimates hold for this global-in-time solution. This finishes the proof of Theorem \ref{Mrt1}.

\section{Conclusion}
In this paper we investigate Rayleigh-Taylor stability problem under Boussinesq approximation, from the viewpoint of rigorous mathematical analysis. As mentioned before, the investigation of such a problem is mainly motivated by recent work on stability problem for MHD system (in the absence of resistivity) under the action of strong background magnetic field. We show that asymptotic stability holds for the motion of incompressible fluid under the action of gravity/bouyancy for the specific steady solution. The main tools used here are spectral analysis (for the linearized operator) and energy method. It seems difficult to extend the method of spectral analysis to other problems such as stability for steady solutions with \emph{variable derivative} temperature, as well as the corresponding problem for nonhomogeneous incompressible model, which will be investigated in our future research.

\centerline{\bf Acknowledgement}
The research is supported by NSFC under grants No. 12071211.


\end{document}